\documentclass[3p,11pt]{elsarticle}
\usepackage{algorithm}
\usepackage{algorithmicx}
\usepackage{algpseudocode}
\usepackage{graphicx}
\usepackage{caption} 
\usepackage{wrapfig}
\usepackage{color,soul}
\usepackage{enumitem}
\usepackage{caption}
\usepackage{subcaption}
\usepackage{comment}
\usepackage{amsmath}
\usepackage{times}
\usepackage{amssymb}
\usepackage{array}
\usepackage{latexsym}
\usepackage{amsthm}
\usepackage{amsfonts}
\usepackage{nomencl}
\usepackage[table]{xcolor}
\usepackage{lineno}

\makenomenclature
\usepackage{hyperref}
\usepackage{array,tabularx}
\usepackage{scalerel,stackengine}
\usepackage{graphics}
\usepackage[english]{babel}
\allowdisplaybreaks
\newtheorem{theorem}{Theorem}[section]

\stackMath
\numberwithin{equation}{section}
\date{}
\newtheorem{l1}{Lemma}[section]

\newtheorem{df}{Definition}[section]

\newtheorem{eg}{Example}[section]
\usepackage[utf8]{inputenc}

\begin{document}

\begin{frontmatter} 

\title{A Radius of Robust Feasibility Approach to Directional Sensors in Uncertain Terrain \tnoteref{mytitlenote}} 
\author[1]{Vanshika Datta}
\author[1]{C. Nahak}
\address[1]{Department of Mathematics, Indian Institute of Technology Kharagpur, Kharagpur, West Bengal, India-721302} 
\cortext[mycorrespondingauthor]{Corresponding author} 
\ead{cnahak@maths.iitkgp.ac.in} 

\begin{abstract}
A sensor has the ability to probe its surroundings; however, uncertainty in its exact location can significantly degrade sensing performance. This motivates the development of coverage strategies that explicitly account for such uncertainty. In this paper, we first incorporate the concept of radius of robust feasibility (RRF) into the coverage modeling of directional sensor networks (DSNs) to systematically handle location uncertainty within a deterministic robust framework. RRF quantifies the maximum allowable perturbation under which feasibility is preserved. We further formalize RRF for DSNs and derive an exact expression tailored to the directional sensing model. Based on this formulation, we then develop a distributed greedy algorithm that optimizes sensor orientations, with the objective of enhancing coverage in the presence of location uncertainty. The proposed approach combines geometric decomposition with robust methodology to ensure reliable performance in the presence of perturbations. Subsequently, we analyze the adaptability of the proposed algorithm in dynamic environments to assess its robustness and efficiency. Moreover, we conduct a series of extensive simulation experiments to evaluate the impact of the proposed methodology and key parameters on the coverage ratio. Experimental results validate its efficacy in maximizing coverage and optimizing sensor orientations, highlighting its potential for practical deployment in uncertain and dynamic environments.

\end{abstract}


\end{frontmatter}

\section{Introduction}
Sensors play a crucial role in modern monitoring systems by enabling the collection of real-time information from physical environments. They are widely used in applications such as environmental monitoring, surveillance, and disaster management, where timely and accurate data is essential \cite{33,35,42}. The performance of such systems largely depends on the quality of sensing, often measured in terms of Quality of Service (QoS), which includes factors like coverage, reliability, and accuracy. Ensuring high QoS is therefore a key objective in the design and deployment of sensor-based systems, with coverage optimization playing a fundamental role, as it directly determines how effectively a region of interest is monitored by the deployed sensors \cite{book}. In practical scenarios, poor coverage leads to coverage holes, redundant sensing, and inefficient resource utilization, which can severely degrade system performance. Higher coverage improves detection capability, reliability, and overall sensing accuracy. Moreover, optimized coverage helps in achieving uniform sensor distribution, reducing energy consumption, and extending network lifetime, while also minimizing deployment costs by avoiding unnecessary redundancy. As a result, coverage optimization has emerged as a central research problem, especially in applications where missing even a small region can lead to significant consequences, such as environmental monitoring, surveillance, and disaster detection \cite{cvg1}. 

Existing works on sensor coverage maximization broadly fall into three categories: Structure-based approaches, Optimization-based methods, and Dynamic and distributed approaches. Structure-based approaches organize the sensing region into well-defined spatial units. These include representative techniques such as Voronoi-based, grid-based, and computational geometry-based methods. Voronoi-based partitioning decomposes the sensing region into local cells, allowing sensors to optimize coverage within their respective regions and enabling scalable implementations \cite{10,39,46}. Grid-based methods discretize the domain into finite points, simplifying evaluation but introducing resolution dependence \cite{30}, while computational geometry approaches provide a precise characterization of coverage through geometric constructs \cite{13}. Optimization-based methods formulate coverage as a mathematical program, whereas metaheuristic approaches such as particle swarm optimization and genetic algorithms explore the solution space efficiently for large-scale and non-convex problems, often prioritizing tractability over optimality \cite{29,36}. Dynamic and distributed strategies include potential field and learning-based methods, where sensors iteratively adjust their positions or orientations based on local interactions. Potential field methods model sensors as particles experiencing attractive and repulsive forces, leading to self-organized deployment, while recent learning-based approaches, such as reinforcement learning, enable adaptive coverage in dynamic environments \cite{15}. 

Wireless sensor networks (WSNs) have emerged as a fundamental technology for large-scale monitoring and have been extensively studied for data acquisition across diverse applications \cite{40}. Most of these works primarily focus on the geometric and algorithmic aspects of sensing. In this context, a majority of existing studies assume omnidirectional sensing models, where sensors provide uniform coverage in all directions, thereby simplifying both modeling and analysis \cite{12}. Under such assumptions, coverage regions are typically represented as disks, allowing tractable geometric and optimization formulations. However, many real-world sensing devices, such as cameras and lidar, inherently operate in a directional manner, making omnidirectional models insufficient for accurately capturing practical scenarios. A more specialized variant, directional sensor networks (DSNs), employs sensors with limited sensing angles, introducing additional complexity in coverage modeling due to orientation-dependent sensing regions \cite{8,43}. Unlike omnidirectional sensing, DSNs require the joint optimization of sensor placement and orientation, making the coverage problem inherently more challenging. Furthermore, coverage overlap, blind spots, and coordination among sensors become more intricate in DSNs compared to omnidirectional settings \cite{37}. This challenge is further amplified in practical deployments, where uncertainties in sensor positioning, calibration errors, and environmental disturbances can significantly degrade coverage performance \cite{24}. In such settings, even small deviations in location or orientation may lead to substantial sensing gaps, highlighting the need for frameworks that explicitly account for uncertainty while preserving effective coverage.

Uncertainty is inherent in WSNs, arising from factors such as imprecise sensor placement, environmental disturbances, and measurement noise, which can significantly affect sensing performance. To address this, a large body of work has relied on probabilistic approaches, where uncertainty is modeled using known distributions and system performance is evaluated in an expected or statistical sense \cite{23}. Probabilistic methods model coverage through detection likelihoods, capturing uncertainty through stochastic representations and expected performance measures \cite{34}. While such methods provide useful insights, they often depend on accurate distributional information, which may not be available or reliable in many real-world scenarios. The perfect-information assumptions in traditional optimization models often lead to suboptimal or infeasible solutions; this highlights the need for methodologies that explicitly address data variability. This limitation has motivated the adoption of robust optimization (RO) techniques, which aim to ensure system performance under worst-case realizations of uncertainty without relying on precise probabilistic assumptions \cite{1}. RO meets this need by incorporating uncertainty explicitly into mathematical models through deterministic uncertainty sets rather than probabilistic distributions \cite{28}. This approach offers a practical, flexible and risk-averse alternative to stochastic methods, particularly valuable when probabilistic data is unavailable or unreliable and in high-stakes, one-time decision-making scenarios. Consequently, RO has emerged as a powerful framework for handling uncertainty, leading to a rich body of literature that traces its development from early foundational models to more advanced formulations \cite{2}.

The foundational idea of RO emerged with Soyster's 1973 work \cite{38}, which introduced a linear framework ensuring feasibility for all data within a convex uncertainty set. This initiated a systematic approach to handling uncertainty in optimization. The field gained significant traction in the late 1990s with the pioneering contributions of Ben-Tal, Nemirovski and El Ghaoui \cite{3,4,17}, who developed computationally tractable methods for both linear and nonlinear optimization problems under uncertainty. Over the years, RO has evolved into a powerful approach for risk-averse decision-making, extending from convex programs to robust discrete optimization, robust mixed-integer and nonlinear formulations and applications in robust linear programming, semidefinite programming and conic-quadratic problems \cite{5,6,25}. Survey studies \cite{16,26} confirmed RO's growing impact across operations research, engineering and finance, with recent advancements expanding into nonconvex and mixed-integer nonlinear programming, broadening its relevance to complex real-world scenarios \cite{32,45}.

RO has addressed distance uncertainty in WSNs \cite{7,41,44}; however its application for DSNs with uncertain sensor locations has remained unexplored. Besides, while RO effectively addresses uncertainty in model parameters by ensuring that solutions remain feasible for all parameter values within a predefined uncertainty set, an important question concerns the extent to which uncertainty can be tolerated without violating feasibility. This leads to the concept of the radius of robust feasibility (RRF), which quantifies the maximum allowable perturbation such that the solution remains feasible \cite{9}. The RRF provides a precise and interpretable measure of robustness, enabling the assessment and comparison of solutions based on their resilience to uncertainty. Initially defined for linear semi-infinite programming by Goberna et al. \cite{19}, it has since been extended to general linear and convex programs, evolving into a powerful analytical tool within RO that links feasibility preservation with structured uncertainty descriptions. 

Though early work focused on ball uncertainty sets, RRF has since been extended to more general structures such as spectrahedra, with key theoretical results on bounds and exact formulations provided by Chuong et al. \cite{11}. RRF is further analyzed for uncertain convex inequalities with multi-convex and compact uncertainty sets, deriving bounds, exact formulas under symmetry conditions, and tractable characterizations for ellipsoids, polytopes, boxes, and the unit ball \cite{18}, while Li and Wang established an exact RRF formula and positivity conditions for general convex programs \cite{27}. The scope of RRF has since expanded to uncertain mixed-integer linear programming, robust convex optimization \cite{20}. A comprehensive survey by Goberna et al. \cite{21} examined its role across continuous and mixed-integer RO. Their work offered a profuse collection of important references for further study.

Recently, Ridolfi et al. \cite{31} demonstrated a significant advancement by directly integrating the RRF concept into optimization modeling, showing its potential for practical decision-making under uncertainty. Despite these advances, existing works on RRF have largely remained confined to abstract feasibility analysis rather than domain-specific modeling, where uncertainty plays a critical role in system performance. In particular, the integration of RRF as a quantitative robustness measure in sensor coverage optimization, where location and orientation uncertainties are inherent, has not been adequately explored. This reveals a clear gap between the theoretical development of RRF and its use in practical, application-driven optimization problems.

Motivated by the above discussion, this article addresses the lack of a unified framework that simultaneously captures directional sensing, location uncertainty, and deterministic robustness guarantees in coverage optimization. This work integrates the concept of RRF within an RO framework to address coverage maximization in DSNs under location uncertainty. Unlike existing approaches that rely on probabilistic assumptions or treat robustness separately, we incorporate the concept of RRF directly into the coverage formulation. In particular, RRF is used to quantify the maximum allowable perturbation in sensor locations while preserving coverage feasibility, and an exact expression is derived for the directional sensing model. This formulation is further embedded within a RO framework to ensure reliable performance under worst-case deviations. To solve the resulting problem, we develop a distributed greedy algorithm that optimizes sensor orientations while accounting for uncertainty, combining geometric decomposition with robustness-driven decision-making. Furthermore, the adaptability of the proposed approach is analyzed in dynamic environments, and extensive simulations are conducted to evaluate its effectiveness under varying conditions and parameter settings. The results demonstrate that the proposed framework significantly improves coverage performance while ensuring robustness to perturbations, highlighting its practical applicability in uncertain and real-world sensing scenarios. Hence, the key contributions of this work are summarized as follows: 

\begin{itemize}
    \item[1] \textbf{Bridging RRF theory to practical DSN coverage:}
    To the best of our knowledge, this is the first work that brings RRF from a predominantly theoretical setting to a practical framework for DSN coverage. It is embedded directly within the coverage model to handle location uncertainty and provide deterministic worst-case guarantees. This represents a significant conceptual shift from abstract robustness analysis to a practical engineering application involving spatial coverage and sensing. 

    \item[2] \textbf{Explicit modeling of sensor location uncertainty:}
    The proposed framework incorporates uncertainty in sensor locations, rather than assuming perfect placement or modeling only sensing uncertainty. By incorporating realistic deployment inaccuracies, it provides a deterministic and practically relevant representation of uncertainty, serving as a complementary alternative to probabilistic approaches when distributional information is limited or uncertain. 

    \item[3] \textbf{Unified DSN–RO formulation:}
    The formulation jointly considers directional sensing, location uncertainty, and robustness. This combination is rare and technically non-trivial. It optimizes sensor orientations in non-convex, orientation-dependent sensing regions, capturing the inherent complexity of DSNs while ensuring reliable worst-case coverage.
    
    \item[4] \textbf{Scalable distributed optimization via Voronoi structure:}
    A Voronoi-based decomposition is combined with a distributed optimization strategy to reduce complexity and enable scalable implementation. This provides an efficient and practical approach for coverage enhancement under uncertainty.
    
    \item[5] \textbf{Comprehensive experimental evaluation:}
    An iterative refinement strategy is used to improve adaptability across varying configurations. Extensive simulations, including parameter studies and comparisons with baseline methods, demonstrate the effectiveness and robustness of the proposed approach.   
\end{itemize}

The remainder of the article is organized as follows: Section 2 presents the necessary preliminaries, including the RRF, Voronoi region construction, the sensing model, and problem description. The RRF for DSNs, forming the basis of the RO framework, is developed in Section 3. Section 4 outlines the proposed model, detailing its design and objectives. The methodology, encompassing the algorithmic framework and implementation, is described in Section 5. Section 6 provides an in-depth experimental analysis to validate the model under various configurations. The concluding remarks with future research directions are discussed in Section 7.

\section{Preliminaries}
This section introduces the formal notation used throughout the paper and establishes the Voronoi-based spatial decomposition, directional sensing geometry, and the adopted sensing model. We begin by presenting the basic notations, definitions, and preliminary results. Let $\mathbb{R}^n$ denote the Euclidean space, $||\cdot||$ the Euclidean norm, and $\mathbb{B}_n$ the open unit ball in $\mathbb{R}^n$. We denote the transpose of a vector $a$ by $a^T$. For $x, y \in \mathbb{R}^n$, the Euclidean distance is $d(x,y)= ||x-y||$. Let $S \subset \mathbb{R} ^n$, later referred to as the set of sensor locations, be a set with at least two elements and a fixed element $ s \in S$. Table~\ref{tab:1} and Table~\ref{tab:2} explain the acronyms and the principal notations used throughout the paper.

\begin{df}(Reference cone) \label{referencecone}
    Given a linear system $\sigma : = \left\{ a_s ^T x \leq b_s, s \in S \right\}$, its reference cone K is defined as: $$K : = cl cone \left\{  \begin{bmatrix}
        a_s\\ b_s
    \end{bmatrix}  , s \in S; \begin{bmatrix}
        0_n \\ 1
    \end{bmatrix} \right\}.$$
\end{df}
Also, note that $\sigma  $ is consistent, if and only if this reference vector does not belong to $K$, i.e.,  $\begin{bmatrix}
    0_n \\ -1
\end{bmatrix} \notin K$.
The system $\sigma $ is called continuous if $S$ is compact and the functions $a: S \to \mathbb{R} ^n, \; a(s) = a_s \text{ and } b: S \to \mathbb{R} , \;b(s) = b_s$, are continuous.\\
A parametric linear system in the face of data uncertainty in its constraints, denoted by $\sigma^{\alpha }  $, can be captured as follows: 
$$\sigma ^{\alpha } := \left\{ a^T_sx \leq b_s , s \in S \right\}; \;(a_i,b_i) \in \mathbb{U} _i ^{\alpha} \subset \mathbb{R} ^{n+1}; \;i=1:p \;,$$
where $(a_i,b_i) $ for $i=1:p$ are uncertain vectors and $\mathbb{U} _i ^{\alpha};\; i=1:p$, which are typically assumed to be compact and convex.\\
The robust counterpart (deterministic form) of the corresponding system $\sigma ^{\alpha} $ with $\mathbb{U_i^{\alpha}}=(\bar{ a_s}, \bar{b_s} ) + \alpha B_{n+1}, \;s \in S, \; i=1:p$ and the corresponding feasible set $F_R^ {\alpha}$ for some $\alpha \in \mathbb{R} $ is:
$$ \sigma _ R ^ {\alpha} := \left\{ a_s^Tx \leq b_s ,\;  \forall (a_s,b_s) \in \mathbb{U_i^{\alpha}},\; s \in S \right\}.$$

\begin{table}[h]
\centering
\caption{List of Acronyms}
\begin{tabular}{ll}
\hline
\textbf{Acronym} & \textbf{Description} \\
\hline
RO & Robust optimization. \\
WSN & Wireless sensor network. \\
DSN & Directional sensor network. \\
RRF & Radius of robust feasibility. \\
ROI & Region of interest (where coverage is evaluated). \\
$VC$ & Voronoi cell (of a particular point from a set of points in 2D). \\
$VD$ & Voronoi diagram (of a set of generator points in $\mathbb{R} ^2$). \\
\hline
\end{tabular}
\label{tab:1}
\end{table}

\begin{table}[h]
\centering
\caption{List of Notations}
\begin{tabular}{ll}
\hline
\textbf{Symbol} & \textbf{Description} \\
\hline
$\mathbb{R}^n$ & $n$-dimensional Euclidean space. \\
$\mathbb{B}_n$ & Unit ball in $\mathbb{R}^n$ with respect to the Euclidean norm. \\
$||\cdot||$ & Euclidean norm. \\
$d(x,y)$ & Euclidean distance between points $x$ and $y$. \\
$m$ & Total number of sensors in the deployment. \\
$\phi_{\omega}$ & Minkowski function associated with a convex set $\omega$. \\
$s$ & A particular sensor for reference. \\
$S = \left\{ s^o_1, \cdots, s^o_m \right\}$ & Nominal location set of sensors. \\
$s^o_i=(x^o_i,y^o_i)$ & Nominal location of the $i^{th}$ sensor in $\mathbb{R} ^2$. \\
$s^w_i$ & Worst case location of  the $i^{th}$ sensor in $\mathbb{R} ^2$. \\
$s^{\rho}_i$ & Robust case location of the $i^{th}$ sensor in $\mathbb{R} ^2$. \\
$o_i$ & Orientation of the $i^{th}$ sensor after ground projection. \\
$\theta_s$ & Angle of view of the sensor. \\
$V_i$ & Selected vertex of orientation of the $i^{th}$ sensor. \\
$A_i(s_i,o_i)$ & Set of all covered points by the $i^{th}$ sensor in $\mathbb{R} ^2$. \\
$r_s$ & Minimum sensing range of a particular sensor $s$. \\
$\mathcal{U}_i^{\alpha}$ & Uncertainty set describing admissible perturbations in $\mathbb{R}^2$. \\
$\mathbb{U}_i^{\alpha}$ & Uncertainty set describing admissible perturbations in $\mathbb{R}^3$.\\
$\alpha$ & Radius of the uncertainty set. \\
$\rho_{i}$ & RRF associated with the $i^{th}$ sensor. \\
$VC(s_i)$ & Voronoi cell of the $i^{th}$ sensor.\\
\hline
\end{tabular}
\label{tab:2}
\end{table}

\subsection{\bf{Construction of Voronoi Regions}}
\begin{figure}[h!]
    \centering
    \includegraphics[width=0.4\textwidth]{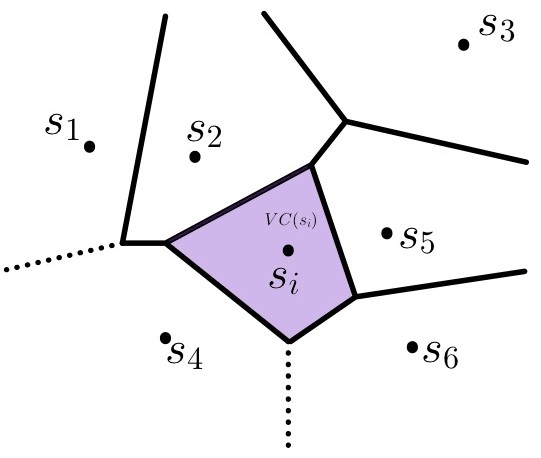} 
    \caption{VD for a set of sensors highlighting $VC(s_i)$}
    \label{fig:voronoidiag}
\end{figure}

A Voronoi diagram (VD) is a spatial decomposition of a region based on proximity to a given set of points, called generators (sensor nodes in our case).
\begin{df}(Voronoi cell) \label{voronoi_cell_def}
    Given a set of generator points $S=\{s^o_1, s^o_2, \dots, s^o_m\}$ in $\mathbb{R}^2$, the Voronoi cell $VC(s^o_i)$ associated with the $i ^{th}$ sensor located at $s^o_i \in S$ is defined as:
    \begin{align*}
        VC(s^o_i) = \{ p \in \mathbb{R}^2 : d(p, s^o_i) \leq d(p, s^o_j) \text{ for all } j \neq i \}, 
    \end{align*}
    where  $j=1:m$, $\;d(a,b)$ is the Euclidean distance between $a$ and $b$.
\end{df}
A VD partitions the domain into regions such that each point belongs to the sensor to which it is closest in Euclidean distance. Thus, every location in the region is uniquely assigned to its nearest sensor node, forming non-overlapping Voronoi cells. Collectively, these cells provide a complete and disjoint coverage of the domain. Accordingly, the VD consists of the union of all Voronoi cells; see Figure~\ref{fig:voronoidiag}.

Geometrically, a Voronoi diagram is constructed by drawing the perpendicular bisectors of the line segments connecting each pair of sensors. These bisectors form the boundaries of Voronoi cells, known as Voronoi edges, while their intersection points are referred to as Voronoi vertices. As a result, with reference of the notations used in Definition \ref{voronoi_cell_def}, the field is partitioned into `$m$' Voronoi cells with the following characteristics:  

\begin{itemize}
    \item[$(1)$] Each sensor $s^o_j$ is contained within a unique cell.
    \item[$(2)$] Any point `$p$' inside the Voronoi cell of the sensor $s^o_j$ is closer to the sensor $s^o_j$ than any other sensor $s^o_t$ for all $j,t \in \left\{ 1,2 , \cdots, m\right\}, \;j \neq t$.
    \item[$(3)$] The Voronoi diagram for a given set of sensors is uniquely determined, meaning that the same set of points will always generate the same Voronoi structure.
\end{itemize}

\subsection{\bf{Radius of Robust Feasibility}}
The RRF $(\rho)$ of an uncertain system $(\sigma ^{\alpha})$ determines the maximum level of perturbation that a system can tolerate while remaining feasible. It is formally defined as follows:
\begin{df}(Radius of robust feasibility) Consider a parametric linear system in face of data uncertainty $\sigma ^{\alpha} $ . Let $F_R^{\alpha}$ denote the feasible set of the robust counterpart of $\sigma ^{\alpha}$. Then, the RRF for the system is given by: 
\begin{align*}
    \rho = \sup \{ \alpha \in \mathbb{R}_+ : (F_R ^ \alpha) \text{ is nonempty} \}.
\end{align*}    
\end{df}
This quantifies the maximum perturbation level under which the system remains robustly feasible. A fundamental concept in computing RRF is the Minkowski function, defined as follows:
\begin{df}(Minkowski function)
    Let $\omega \subset \mathbb{R} ^ n $ be a convex set containing $0_n$ in its interior. Then the Minkowski or gauge function of $\omega$ denoted by $\phi_{\omega}$, where  $\phi_{\omega} : \mathbb{R} ^n \to \mathbb{R} _+ := [0, + \infty [ $ is given by: 
    $$\phi_{\omega} (x) := \inf{ \left\{ t>0 : x \in t \omega \right\} },\; x \in \mathbb{R} ^n.$$
\end{df}
The following lemma provides some properties of the Minkowski function.
\begin{l1}[\cite{11}, Lemma 1.3.13] Let $\omega \subset \mathbb{R} ^n $ be a convex set such that its interior contain $0_n$, then, the following properties hold:
\begin{itemize}
    \item[$(1)$] $\phi_{ \omega}$ is sublinear and continuous.
    \item[$(2)$] $\left\{ x \in \mathbb{R} ^n : \phi _{\omega} (x) \leq 1 \right\} = cl (\omega)$, where $cl (\omega ) $ stands for the closure of $\omega$.
    \item[$(3)$] If in addition, $\omega$ is bounded and symmetric, then, $\phi _ {\omega }: = || \cdot ||$ is a norm on $\mathbb{R} ^n$ generated by $\omega$.
\qed
\end{itemize} 
\end{l1}

Exact formula of RRF for an arbitrary set $S$ has been obtained in literature in terms of the so-called hypographical set $H(\bar a, \bar b)$ of the nominal system $\sigma$: 
$$H(\bar a, \bar b) := conv \left\{ (- \bar{a_i} , -\bar{b_i}): i=1:p \right\} + \mathbb{R} _+ (0_n, -1),$$
where $\bar a = ( \bar{a_1},\bar{a_2}, \cdots, \bar{a_p} ) \in (\mathbb{R} ^n)^p \text{ and } \bar b = ( \bar{b_1},\bar{b_2}, \cdots, \bar{b_p} ) \in \mathbb{R} ^p$.
\begin{theorem} \cite{11}
    If the nominal system is feasible, the RRF of $\sigma^{\alpha}$ is:
    \begin{align*}
        \rho = \inf_{(a,b) \in H(\bar{a},\bar{b})} \phi_Z (a,b),
    \end{align*}
    where $Z \subset \mathbb{R} ^n $ is a convex and compact set containing $0_n$ in its interior. 
\end{theorem}

\subsection{\bf{Sensing Model in Directional Sensor Networks}}\label{2.3}
\begin{figure}[h!]
    \centering
    \includegraphics[width=0.5\textwidth]{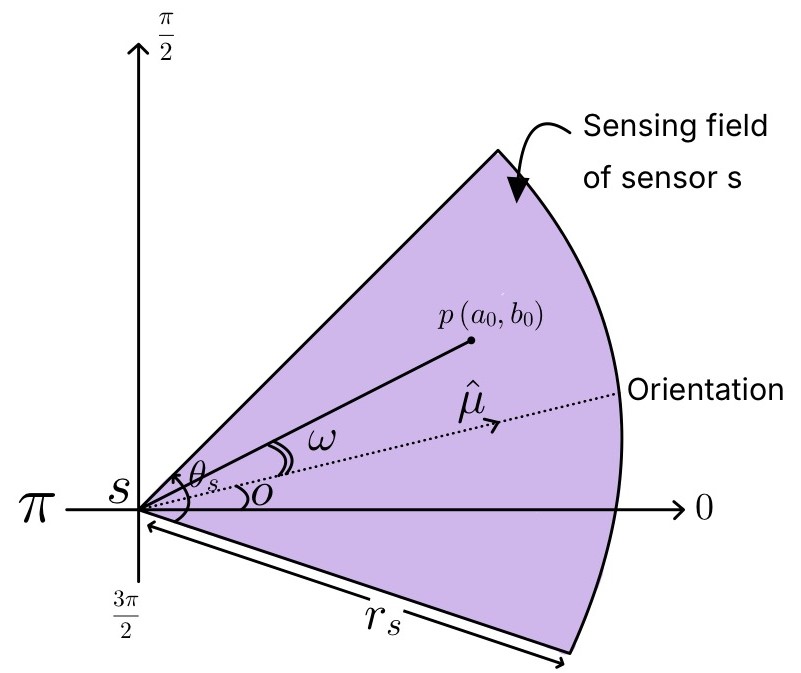}
    \caption{Sensing model for directional sensors}
    \label{fig:s_model}
\end{figure}

In DSNs, each sensor has a limited sensing angle and range, as shown in Figure~\ref{fig:s_model}.
Here, $r_s$ is the sensing range of the sensor $s$, located at $(x_s,y_s)$. $\hat{\mu }$ is a unit vector denoting the orientation of the sensor. Note that each sensor has only one orientation, making an angle '$o$' lying between $-\pi$ and $ \pi$ relative to the positive $x$ direction. The effective sensing field, called the field of view of the sensor, is a sector denoted by the shaded part in the figure. It is formed using two side lengths of sensing sector with radius $r_s$ and an included angle $\theta _s$, which is the angle of view or effective viewing angle.\\
Note that the sensor `$s$' positioned at $(x_s, y_s)$ with sensing range $r_s$ and sensing angle $\theta_s$ covers a point $p(a_0, b_0)$ if the following two conditions are satisfied:
\begin{itemize} 
    \item[$(1)$] The Euclidean distance between point $p$ and sensor $s$ is less than or equal to the sensing radius $r_s$ of sensor `$s$', i.e.,
    $$\sqrt{(x_s-a_0)^2+(y_s-b_0)^2} \leq r_s.$$
    \item[$(2)$] The included angle $\omega $ between $\vec{sp}$ and $\hat{\mu}$ is less than or equal to half of the angle of view $\theta _s$ of sensor `$s$', i.e.,
    $$\omega \leq \frac{\theta _{s}}{2}. $$
\end{itemize}

\subsection{\bf{Problem Formulation and Solution Approach}}
In large-scale deployments, sensors are typically placed in an ad hoc or random manner over extensive regions, resulting in location data that is only known approximately rather than precisely. Such positional uncertainty is often either ignored or modeled probabilistically in the literature, which may lead to coverage solutions without guaranteed feasibility under all realizations of uncertainty. 

In this work, we address the problem of maximizing coverage while minimizing cost, i.e., reducing the number of sensors, subject to a guaranteed minimum coverage level even in the presence of location uncertainty. To achieve this, we incorporate the RRF into the DSN framework to explicitly account for bounded uncertainty in sensor positions. The RRF quantifies the maximum allowable deviation of each sensor from its nominal location while preserving feasible coverage. Building on this, we formulate an RO problem that maximizes worst-case coverage by adaptively adjusting sensor orientations. To solve this problem efficiently, we propose a distributed greedy algorithm that iteratively refines sensor orientations under uncertainty, ensuring high coverage performance with minimal resource utilization.

\section{Radius of Robust Feasibility for Directional Sensor Networks}
In this section, we extend the RRF framework to DSNs by grounding uncertainty directly in sensor misplacement. 

Instead of arbitrary coefficient perturbations, we restrict changes in $(a_s, b_s)$ to those induced by sensor location errors. Specifically, the perturbations in the system $\sigma$ (see Definition \ref{referencecone}) arise from the uncertainty in sensor nodes $s \in S$. Thus, we consider only those variations in $(a_s, b_s)$ in $\sigma^\alpha$ that result from uncertain sensor locations. In this way, we set up the following notation.\\
For a particular sensor location in $S = \left\{ s_1^o,s_2^o, \cdots, s_m^o \right\} $, let the nominal location be $s_i^o$ and let the uncertainty set capturing all feasible perturbations be $\mathcal{U} _i ^{\alpha} \;; \; i=1:m$, defined by: 
$$\mathcal{U} _i ^{\alpha} := {s_i^o} + \alpha \mathbb{B} _2 ; \; i=1:n\;, $$
where $\mathbb{B}_2$ denotes the closed unit ball in $\mathbb{R} ^2$, which is convex and compact with $0_{2} \in \operatorname{int}(\mathbb{B}_2)$. This structured formulation allows the derivation of robustness conditions that are physically grounded in sensor placement errors. Assuming the nominal system $\sigma$ is feasible, we now define the uncertain sensor position `$s_i$' for sensor located at $s_i^o \in S$ along the direction of sensing $\vec{u}$, which gives a meaningful and tractable representation for robust coverage analysis: $$s_i = {s_i^o} + \alpha \frac{\vec u}{|\vec u|}.$$

Then, for a fixed realization $s_i = ({x_i}, {y_i}) \in S$, the Voronoi cell VC$(s_i)$, by Definition \ref{voronoi_cell_def}, corresponds to the system of inequalities:$$VC(s_i) := \left\{ x \in \mathbb{R} ^2 : ||x-s_i||_2^2\leq ||x-s_j||_2^2 \;, \;\forall \; j \neq i \; \right\}.$$
Expanding the norms and rearranging terms, the inequality becomes $$ VC(s_i) = \left\{ x \in \mathbb{R} ^2 : 2 (s_j-s_i)^T x \leq ||s_j||_2^2 - ||s_i||_2^2,\; \forall \; j \neq i \right\} ,$$ which defines the half-space constraints forming the Voronoi cell of $s_i$.
We aim to express the induced uncertainty in the coefficients $(a_{ij}, b_{ij} )$ due to the uncertainty in $s_i$. Let us define: 
\begin{equation} \label{a_and_b}
    a_{ij}:= 2(s_j-s_i) \in \mathbb{R} ^2 \text{ and } b_{ij} := ||s_j||_2^2 - ||s_i||_2^2 \in \mathbb{R}.
\end{equation}
Then $a_{ij}^Tx \leq b_{ij}$ becomes the Voronoi half-space system of constraints under uncertainty. This leads to an uncertainty set in the coefficients $(a_{ij}, b_{ij}) $, which are functions of $(s_i,s_j)$.
Next, we consider the set of all such uncertain locations which are contained in the Voronoi cell of $s_i$. Note that the feasible set of such balls is the same as the following feasible set: 
\begin{equation}\label{robust_feasible_set} 
\mathcal{F} ^{\alpha} = \left\{ x \in \mathbb{R} ^2 :  a_{ij}(s_i,s_j)^Tx\leq b_{ij}(s_i,s_j) \;, \forall \;j \neq i, \forall s_i \in \mathcal{U} _i ^{ \alpha} , \forall s_j \in \mathcal{U} _j ^{ \alpha}   \right\}.
\end{equation}
Define the uncertainty set for the coefficients as:
\begin{equation}\label{nonlinear}
\begin{aligned}
    \mathbb{U} _{ij} ^ {\alpha} : =  \big\{ (a_{ij},b_{ij}) \in \mathbb{R} ^3: \; &  a_{ij}:= 2(s_j-s_i) \in \mathbb{R} ^2 , \;\\& b_{ij} := ||s_j||_2^2 - ||s_i||_2^2, \forall s_i \in \mathcal{U} _i ^{ \alpha} , \forall s_j \in \mathcal{U} _j ^{ \alpha}   \big\}.
\end{aligned}    
\end{equation}
This defines a nonlinear uncertainty set in $(a,b)$, induced from the ball uncertainty in $s_i$. The nonlinear uncertainty set in $(a,b)$ must be approximated linearly for tractability and compatibility with RO frameworks. Hence, to align our model with the linear uncertainty framework in the literature, we use the Taylor expansions to approximate the nonlinear functions of $s_i$ and $s_j$ that define $a_{ij}$ and $b_{ij}$, thereby obtaining linear expressions that can be used to define linear uncertainty sets. 
To this end, we assume $s_i = {s_i^o} + \Delta s_i$ and $s_j = {s_j^o} + \Delta s_j, \text{ with } || \Delta s_i || , || \Delta s_j || \leq \alpha$. Also, let us define using equation~(\ref{a_and_b}) that $a= 2(s_j-s_i)$ and $ b = ||s_j||^2-||s_i||^2$ have their nominal values as: ${a^o}= 2({s_j^o}-{s_i^o})$ and $ {b^o} = ||{s_j^o}||^2-||{s_i^o}||^2$. Hence, we have $$a = {a^o} + 2(\Delta s_j - \Delta s_i). $$
Rearranging this, we obtain $$||a-{a^o}|| \leq 2 (|| \Delta s_j|| + || \Delta s_i ||) \leq 4 \alpha .$$
Similarly, we have: $b = ||s_j||^2-||s_i||^2 $. Utilizing first-order Taylor expansion, we get $$||s_j||^2 \approx ||{s_j^o} ||^2 + 2 {s_j^o}^T \Delta s_j,$$ and $$||s_i||^2 \approx ||{s_i^o} ||^2 + 2 {s_i^o}^T \Delta s_i.$$
This implies that $$b \approx {b^o} + 2 ({s_j^o}^T \Delta s_j-{s_i^o}^T \Delta s_i).$$
Now, we bound the deviation from ${b^o}$: 
\begin{align*}
    ||b - {b^o}|| & \leq 2 (||{s_j^o}^T \Delta s_j||+||{s_i^o}^T \Delta s_i||)\\
    & \leq 2(||{s_j^o}|| \cdot || \Delta s_j|| + || {s_i^o}|| \cdot || \Delta s_i ||) \\
    & \leq 2 \alpha (||{s_j^o}||+ ||{s_i^o}||)\\
    & := \delta_b(\alpha).
\end{align*}
Further, we define the linear uncertainty set for a particular $i\in \left\{ 1, 2, \cdots, n \right\}$ using equation~(\ref{nonlinear}) as follows: 
\begin{equation} \label{linear}
    \mathbb{U} _i = \left\{ (a,b) : || a - {a^o} || \leq 4 \alpha,\; ||b - {b^o}| | \leq \delta_b(\alpha) \right\}, \forall j=1:n, j \neq i.
\end{equation}
Hence, the robust feasible region from equation~(\ref{robust_feasible_set}) and equation~(\ref{linear}) becomes $$ \mathcal{F}^{\alpha} = \left\{ x \in \mathbb{R} ^2 |  a ^T x \leq b, \; \forall (a,b) \in \mathbb{U} _ i\; \forall j \neq i \right\} ,$$ yielding a standard tractable uncertainty model, compatible with RO and we can compute the RRF for this using similar approach as in \cite{11}. Hence, we define the robust feasible region for a sensor located at $s_i \in S $ as the intersection of these uncertain half-spaces. Moreover, we define the notion of RRF for DSN as follows:
\begin{df}(RRF for DSN) 
    Let $S$  be a set of sensor locations in a DSN. The radius of robust feasibility for directional sensor located at $s_i$ is the largest uncertainty radius $\alpha$ such that its robust feasible region $\mathcal{F}^{\alpha}$ remains nonempty, i.e., $$\rho_i = \sup{ \left\{ \alpha \in \mathbb{R} \geq0 : \mathcal{F}  ^{\alpha}  \neq \phi  \right\} }\; ;\; \text{for } i=1:m.$$
\end{df}
The RRF defines the largest uncertainty radius ensuring feasibility. In sensor networks, it accounts for minor displacements from nominal positions and guides the maximum permissible deviation when exact deployment is infeasible.\\
Hence, we state the first main result of the article in the following theorem, which gives an exact formula for the RRF of a sensor in a DSN. 
\begin{theorem}(Exact RRF formula) \label{theorem_main}
Let the nominal system $\sigma$ be feasible. Then, the RRF of the $i^{th}$ sensor `$s_i$' in a DSN is: $$\rho_i = \inf_{x \in \mathcal{F} ^{\alpha}} \min_i \frac{{b_i^o}- {a_i^o}^T x }{\xi_{\mathbf{Z} _i} (x,-1)}\;;\; i=1:m$$
 where $\xi_{\mathbf{Z} _i} (x,-1)$ is the support function of the uncertainty set $\mathbf{Z} _i$, capturing worst-case directional perturbation.
 \end{theorem}
 \begin{proof}
 We begin with the definition of RRF for a system with uncertain linear constraints: $$(a_i,b_i) \in ({a}_i^o,{b}_i^o) + \rho_i  \mathbf{Z} _i \;;\; i = 1:m .$$ 
    We require the system to remain feasible under all such perturbations, i.e., to find $\rho_i$ such that
    $$ \bigcap_{(a_i,b_i) \in ({a}_i^o,{b}_i^o) + \rho_i \mathbf{Z}_i} \{ x \mid a_i^T x \leq b_i \} \neq \emptyset.$$
    This is the same as $$\forall i, \; \sup_{(a_i,b_i) \in ({a}_i^o,{b}_i^o) + \rho_i \mathbf{Z} _i } (a_i^T x - b_i) = \rho \cdot \xi_{\mathbf{Z} _i } (x,-1) \leq {a_i^o}^T x - {b_i^o} , $$
    where $\xi_{\mathbf{Z}_i}(x, -1) $ denotes the support function of the uncertainty set $ \mathbf{Z}_i $, evaluated in the direction $(x, -1)$.
    We have $$\rho_i \leq  \frac{{b_i^o}- {a_i^o}^T x }{\xi_{\mathbf{Z} _i} (x,-1)} .$$
    To ensure all constraints remain satisfied under uncertainty, we consider the minimum over all $i$ and then determine the maximum feasible radius by taking the infimum over all feasible $x \in \mathcal{F}  ^{\alpha}$. As a result, we deduce $$\rho_i =  \inf_{x \in \mathcal{F} ^{\alpha}} \min_i \frac{{b_i^o}- {a_i^o}^T x }{\xi_{\mathbf{Z} _i} (x,-1)},$$ 
    which is the required expression, thereby concluding the proof.
\end{proof}
Note that the expression in Theorem~\ref{theorem_main} quantifies the robustness margin $\rho_i$ for the $i^{th}$ sensor $s_i$ by minimizing over positions $x$ within its nominal Voronoi cell while accounting for the deviations in each constraint. The formula enables tractable computation of robust Voronoi regions even under spatial uncertainty.
Next, we state the following theorem, which will be helpful in our application of DSN area coverage for selecting the working direction of a particular sensor.
\begin{theorem} \label{t4}
The working direction of a sensor that maximizes coverage within its corresponding Voronoi cell is oriented toward a vertex of that Voronoi cell.
\end{theorem}
\begin{proof}
Consider a Voronoi diagram $VD(S)$ generated by a finite set of sensor locations $S$ in a Euclidean plane. Each Voronoi cell is a convex polygon whose vertices are equidistant from at least three sensors.\\
Let us assume, if possible, there exists a point $P$ in a Voronoi cell that is farther from its associated sensor than any vertex of the cell. Since $P$ is not a vertex and is farther than a vertex, it must lie on the edge of the Voronoi cell. However, this assumption leads to a contradiction based on the following reasoning:
\begin{enumerate}
    \item[$(1)$] Due to the convexity of the Voronoi cell, any infinitesimal displacement of $P$ along the edge, either left or right, remains within the cell. Consequently, no such perturbed point can exceed $P$ in distance from the sensor, contradicting the assumption that $P$ is the farthest point.
    \item[$(2)$] If both adjacent points to $P$ on an edge are equidistant from the sensor, the edge would form a circular arc centered at the sensor. However, Voronoi edges in a Euclidean diagram with finitely many sensors are always straight-line segments, leading to a contradiction. Thus, such a point $P$ cannot exist.
\end{enumerate}
Hence, the farthest point within a Voronoi cell must always be at a vertex, proving that the optimal working direction of a sensor is towards a vertex.
\end{proof}

\section{Proposed Model}
In this section, we begin by describing our proposed sensing model, then formulate the corresponding optimization problem, outline the key constraints and stopping criteria, and conclude with a discussion of the main assumptions and structural properties of the problem.

\subsection{\bf{Sensing Model}}
Sensor locations may deviate from their nominal points. For each sensor, we compute the maximum uncertainty radius that preserves feasibility and orient the sensor accordingly. Only the portion inside its Voronoi cell is considered to reduce coordination complexity and avoid redundant overlaps. Although sensor locations are uncertain within an RRF neighborhood, Voronoi cells constructed from nominal positions provide a stable responsibility region for each sensor. This ensures localized optimization and avoids global conflicts while still accounting for deviations. Our model primarily consists of a DSN, where each sensor has specific coverage areas defined by its orientation and effective range. For a particular sensor `$s$', the effective coverage area is given by $$A_s = \theta_s \cdot r_s^2,$$ where $\theta_s$ is the angle of coverage in radians and $r_s$ is the effective coverage range. We consider the intersection of this coverage area with the Voronoi region of `$s$'.

\begin{figure}
\centering
\begin{subfigure}{0.48\columnwidth}
    \centering
    \includegraphics[width=0.8\linewidth]{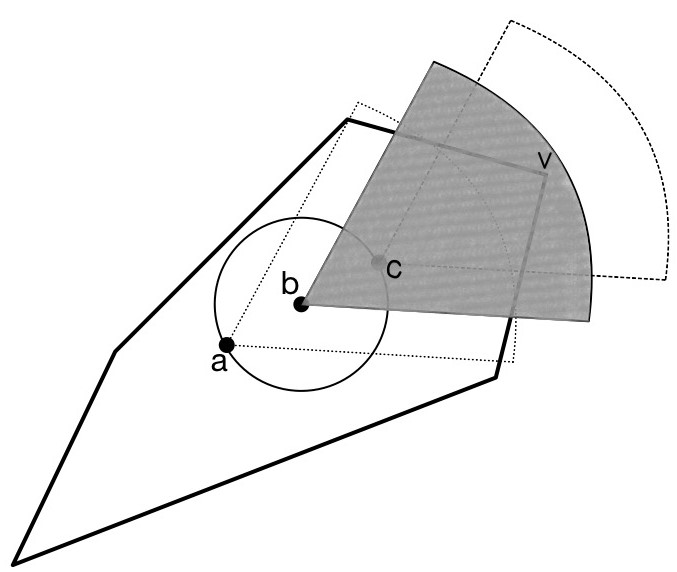}
    \caption{Best and worst-case coverage of a sensor}
\end{subfigure}
\hfill
\begin{subfigure}{0.48\columnwidth}
    \centering
    \includegraphics[width=\linewidth]{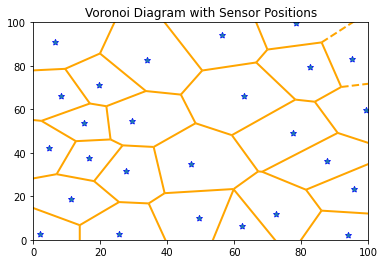}
    \caption{Construction of Voronoi diagram}
\end{subfigure}
\caption{Illustration of robust coverage and Voronoi construction}
\label{fig:best_worst_in_cell}
\end{figure}

Figure~\ref{fig:best_worst_in_cell}a illustrates the potential best and worst-case sensor placement scenarios when facing a vertex `$v$'. The uncertainty in sensor positioning is represented by a circular region centered at its nominal location `$b$', within which the sensor may be positioned. This variability affects the effective coverage area. Each sensor may lie anywhere in its uncertainty disc; the approximated worst-case position is the one minimizing cell-wise coverage, and the orientation is chosen toward the corresponding critical Voronoi vertex. To ensure robust performance, we assume the most disadvantageous sensor placement while determining its optimal orientation. The shaded region denotes the nominal area it covers; the bounded dotted and dashed regions denote the area it covers under the potential best and worst-case scenarios, respectively. Thus, the sensor is considered to be at location `$c$', where its coverage within its Voronoi cell is minimized. If the optimal orientation aligns towards vertex `$v$', the conservative guaranteed coverage will be at least that achieved when the sensor is positioned at `$c$', while in the potential best-case scenario, it could be maximized if the sensor is located at `$a$'. Finally, the coverage is obtained by averaging over all vertex-wise potential worst-case placements across the Voronoi cell, providing a robust estimate for each sensor.

\begin{eg}[Single sensor with positional shift]
Let a sensor `$s$' have a coverage range of $s_r=100 $ m and a positional uncertainty of $\rho_s = 5$ m. The robust coverage range is then given by $r_s = s_r - \rho_s$, resulting in 95 m. Since this represents the potential worst-case scenario, the sensor can ensure coverage of at least 95 m within its Voronoi cell when oriented toward its chosen vertex.
\end{eg}

\subsection{\bf{Optimization Formulation}} 
The original robust constraints are semi-infinite since coverage feasibility must hold for all points in a continuous region. We reduce these to a finite set of critical constraints using Voronoi vertices, boundary intersections, and sector arc points, which correspond to the extremal directions for robust coverage feasibility. Given a bounded region $R$ and $m$ sensors, we construct a Voronoi diagram (Fig.~\ref{fig:best_worst_in_cell}(b)) where each sensor owns a convex cell, and each point in the cell is closest to that sensor. 

Based on the characteristics of a Voronoi cell, the sensor within a given cell has the strongest inductive capacity or the best sensing quality for any position inside that cell. Therefore, each sensor should be responsible for covering its respective cells as fully as possible. Additionally, since Voronoi cells are always convex polygons, Theorem~\ref{t4} states that the farthest point within a cell from its associated sensor is typically one of the cell's vertices.

Since we focus on perturbations caused by sensor location uncertainty rather than arbitrary system coefficient changes, we consider only variations arising when the sensor positions $s_i \in S\; ; \; i =1:m$ are unknown and subject to uncertainty. Thus, we introduce the following notations.

Let $I_i$ be the index set of Voronoi vertices for the $i^{th}$ sensor located at $s_i$ for $ i =1:m$, where $n_i$ denotes the number of vertices in the Voronoi cell of $s_i$. The vertices of the Voronoi region associated with the $i^{th}$ sensor as $v_{ij}$ for $j=1:n_i$. The area covered by $s_i$ within its Voronoi region when oriented towards vertex $v_{ij}$ is represented as $A_{ij}$. Furthermore, let $o_i$ and $A_i(s_i,o_i)$ denote the optimal orientation selected and the corresponding covered area for the sensor located at $s_i$ within its Voronoi region, respectively, for $i=1:m$. 

To maximize the total coverage while ensuring robustness against positional uncertainty, we formulate the following optimization problem:
\begin{align*}
    \max_{o_i} \sum_{i=1}^{m} A_i(s_i,o_i),
\end{align*}
subject to constraints ensuring that the sensor orientations account for the approximated worst-case positioning within the uncertainty set, thereby guaranteeing robust feasible coverage across the entire network.

Hence, we have the following constraints:
\begin{itemize}
    \item[$(1)$] {\bf Initial Orientation Selection:} Each sensor selects the vertex that provides the maximum coverage within its Voronoi cell:
    \begin{align*}
        o_i = \arg\max_{v_{ij} \in I_i} A_{ij}, \quad \forall j \in \{1,2,\dots,n_i\}.
    \end{align*}
    \item[$(2)$] {\bf Overlap Resolution Constraint:} If two sensors select the same vertex, the sensor with the lower coverage area selects the next best option:
    \begin{align*}
        o_i' = \arg\max_{\substack{v_{ij} \in I_i \\ v_{ij} \neq o_i}} A_{ii}, \quad \text{if } o_i = o_k \text{ and } A_i(s_i,o_i) < A_k(s_k,o_k).
    \end{align*}
    \item[$(3)$] {\bf Boundary Avoidance Constraint:} Sensors should not orient towards vertices that are within an $\epsilon$-distance of the boundary $\partial R$:
    \begin{align*}
        o_i \notin \{ v_{ij} \in I_i \mid d(v_{ij}, \partial R) < \epsilon \}, \quad \forall i=1:m.
    \end{align*}
    \item[$(4)$] {\bf Uncertainty Constraint:} Each sensor's actual location $(x_i, y_i)$ lies within an uncertainty ball of radius $r_i$ around its nominal position $s_i^o$:
    \begin{align*}
        (x_i, y_i) \in B(s_i^o, r_i), \quad \forall i=1:m.
    \end{align*}
\end{itemize}
where $(x_i,y_i)$ denotes the actual position of the $i^{th}$ sensor, constrained within its uncertainty set. 
This optimization model ensures that:
\begin{itemize}
    \item[$(1)$] Sensors are oriented for potential maximum coverage within their Voronoi regions.
    \item[$(2)$] Overlapping orientations are resolved dynamically by reassigning the sensor with lower coverage.
    \item[$(3)$] Sensors avoid aligning towards vertices near the boundary for better stability.
    \item[$(4)$] The model accounts for uncertainty in sensor positioning, ensuring robust feasible coverage.
\end{itemize}

\begin{figure}[t]
\centering
\begin{subfigure}{0.48\columnwidth}
    \centering
    \includegraphics[width=\linewidth]{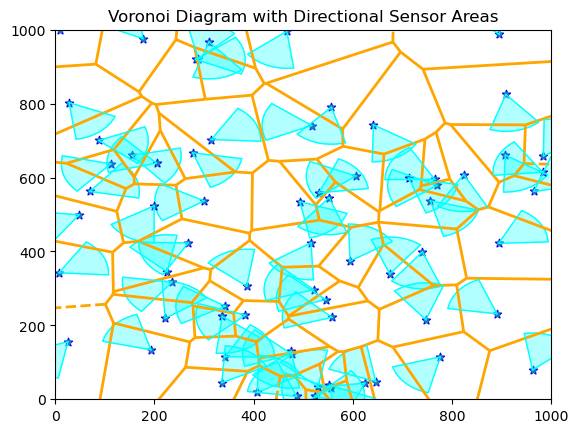}
    \caption{Nominal Sensors}
\end{subfigure}
\hfill
\begin{subfigure}{0.48\columnwidth}
    \centering
    \includegraphics[width=\linewidth]{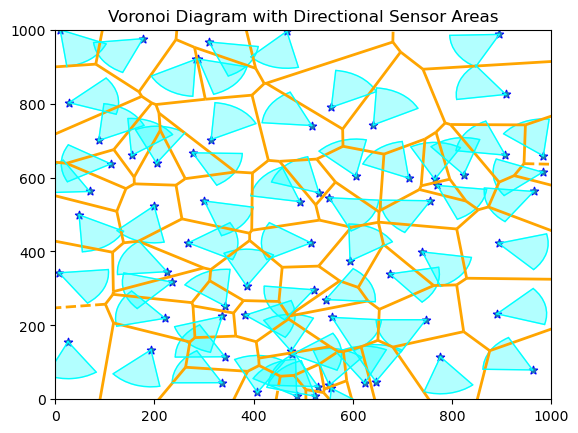}
    \caption{Robustified Sensors}
\end{subfigure}
\caption{Comparison of sensor orientation strategies}
\label{fig:orientation_comparison}
\end{figure}

If a sensor at location $s_i$ changes its orientation due to overlap resolution and all vertex options are exhausted, the sensor can either be put into sleep mode to conserve energy or select the initially preferred vertex. The following example shows that even if the sensor selects the initially preferred vertex, it can still contribute to the robustness of the overall model. As we shall see in the following sections, Figure~\ref{fig:orientation_comparison} illustrates the contrast between nominal and robustified orientations under this integrated adjustment process, highlighting the effect of uncertainty-aware optimization on coverage stability. It can be observed that the robustified orientations provide improved coverage stability by accounting for positional uncertainty, even in the presence of overlap.

\begin{eg}[Overlapping Sensors with RRF Constraint]
Suppose sensors $s_1$ and $s_2$ are placed 150 m apart and both oriented toward the same vertex $V$, with robust ranges of $r_{s_1} = 80$ m and $r_{s_2} = 90$ m. This results in an overlapping coverage area, leading to unnecessary redundancy. To mitigate this, the sensor contributing lesser effective coverage within its Voronoi cell, say $s_1$, could be reoriented toward an alternative vertex. However, if all other vertices are already assigned to sensors with coverage greater than that of the selected sensor, then $s_1$ may still remain oriented toward $V$, ensuring robustness against positional uncertainty. The overlap, $\epsilon = 80+90-150 =20$ m, acts as a buffer against deployment errors, preventing coverage gaps and maintaining reliable sensing.\\
{\bf Visual representation:} Imagine two sectors representing the feasible regions of $s_1$ and $s_2$, slightly overlapping. Even if $s_1$ or $s_2$  moves slightly within their feasible radii, the overlapping area remains covered, ensuring robust coverage continuity.
\end{eg}

A connectivity or energy constraint can also be incorporated into the constraint set; however, we omit it here as our focus is solely on maximizing coverage while accounting for location uncertainty. Additionally, a constraint ensuring that the covered area exceeds a threshold $\epsilon$ could be introduced; otherwise, the sensor can be put into sleep mode to conserve energy.

\subsection{\bf{Geometric Considerations and Practical Assumptions}} 
To streamline our analysis, we assume the region is bounded, with positional uncertainty modeled as a ball centered at each sensor's nominal location. Additionally, as discussed earlier, we focus exclusively on perturbations arising from uncertain sensor positional shifts, rather than general coefficient variations.

\begin{figure}
    \centering
    \includegraphics[width=0.4\textwidth]{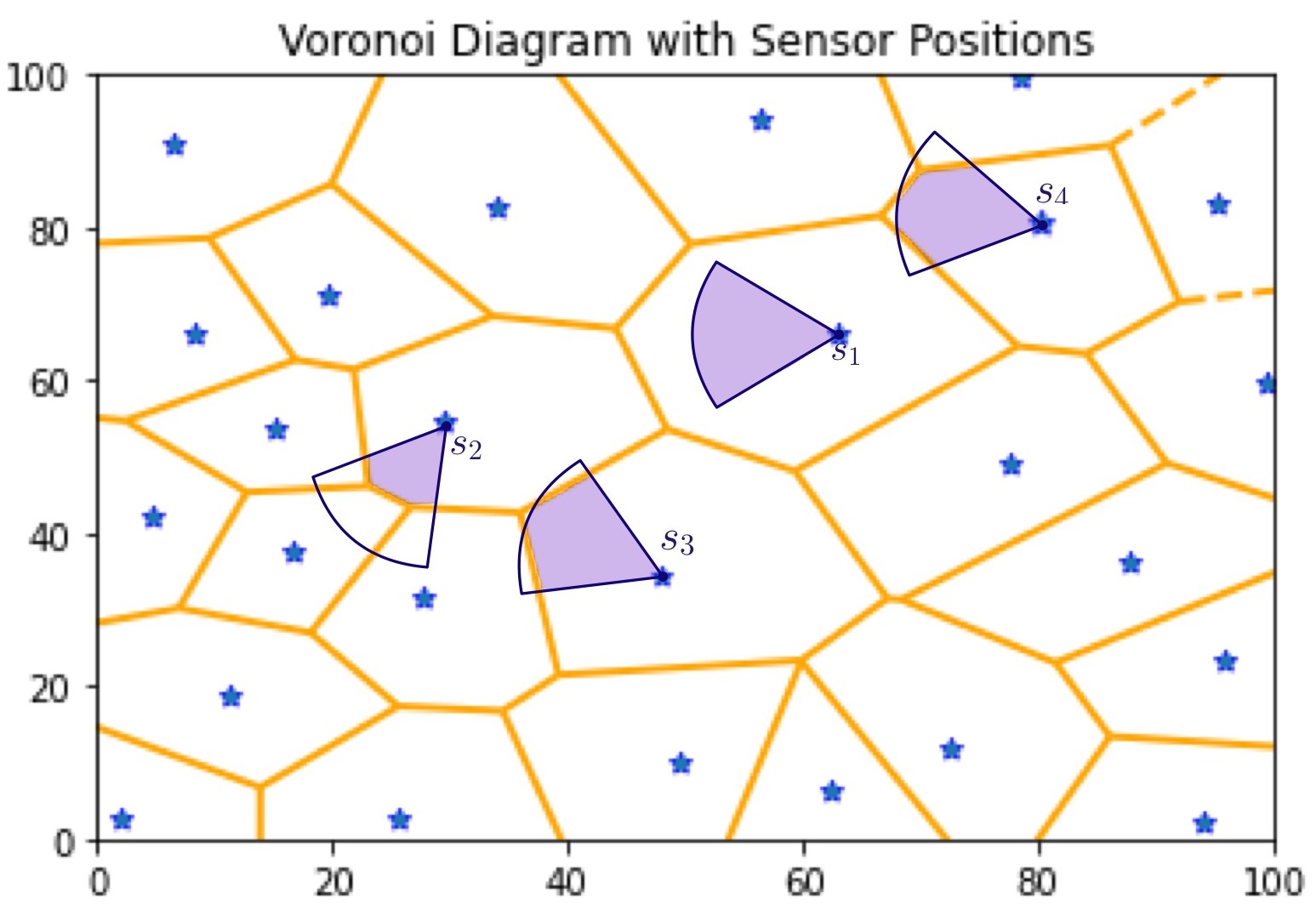}
    \caption{Possible cases for area coverage}
    \label{fig:geometric}
\end{figure}

In Figure~\ref{fig:geometric}, sensors exhibit different coverage patterns within their respective Voronoi cells. As Voronoi cells are convex, coverage patterns reduce to four distinct cases, expressed in Figure~\ref{fig:voronoi_diag}, into which all other configurations can be categorized. 

\section{Proposed Methodology}
We utilize Voronoi diagrams and direction-adjustable sensors to design an adaptive, distributed placement strategy that improves coverage in DSNs. As outlined in Section~\ref{2.3}, a point is covered if it satisfies two conditions, which are evaluated using the sensor’s robust location rather than its nominal ones. Consider the $i^{th}$ sensor $s_i;\; i = 1:m,$ be located at the nominal coordinates $(x_i^o,y_i^o)$, with an RRF $\rho_i,$ and oriented toward vertex $V=(a,b)$. The direction vector from $s_i$ to $V$ is given by $\vec{\mu} = (a-x_i,b-y_i)$. After normalization and scaling by the RRF, the robust location of sensor $s_i^o$ is given by $s_i^w=(x_i^w,y_i^w)$ such that
\begin{align}
    (x_i^w,y_i^w) &= (x_i^o,y_i^o) + \rho_i \frac{\vec{\mu}}{|\vec{\mu}|} \notag \\
    &= (x_i^o,y_i^o) + \rho_i \left(\frac{a-x_i^o}{\sqrt{(a-x_i^o)^2 + (b-y_i^o)^2}},\frac{b-y_i^o}{\sqrt{(a-x_i^o)^2 + (b-y_i^o)^2}}\right) \label{5.2}.
\end{align}
Thus, a point $p(a_0,b_0)$ is considered covered by a particular sensor located at $s_i^o$ if it meets the following conditions: $d(s_i^w, p) \leq r_i$ and $\omega \leq \frac{\theta_s}{2}$, where $r_i$ and $\theta_s$ denote the effective robust coverage range and effective angular coverage range of the $i^{th}$ sensor, respectively. The parameter $\omega $ represents the angle between the vector $\vec{sp}$ and the unit vector $\hat{\mu}$.\\
Note that if $o_i $ is the angle between $\hat{\mu} $ and the relative positive x-axis with $-\pi \leq o_i \leq \pi$, one can derive the following from the second condition of Section~\ref{2.3}:
\begin{align}
    \vec{sp} \cdot \hat{\mu} &= ||\vec{sp} ||\; ||\hat{\mu}|| \; \cos{\omega} \notag\\
    \implies (a_o-x_i^w) \cos{o_i}+(b_0-y_i^w) \sin{o_i}&=\sqrt{(a_0-x_i^w)^2 + (b_0-y_i^w)^2} \cos{\omega} \notag \\ \label{cond2}  \implies (a_o-x_i^w) \cos{o_i}+(b_0-y_i^w) \sin{o_i}&\geq \sqrt{(a_0-x_i^w)^2 + (b_0-y_i^w)^2} \cos{\frac{\theta_s}{2}}. 
\end{align}
Thus, a point $p(a_0,b_0)$ is considered covered by a particular sensor at $s_i^o$ under robust case if it satisfies inequalities $d(s_i^w, p) \leq r_i$ and equation~(\ref{cond2}). Note that the proposed formulation optimizes coverage within each sensor’s Voronoi cell, as this provides a natural spatial decomposition for local maximization. However, in computing the total network coverage, we explicitly account for the fact that points within a given Voronoi cell may also be covered by neighboring sensors, thereby ensuring that overlapping coverage across cells is properly incorporated.

\subsection{\bf{Localized Orientation Optimization}} 
Each sensor initially determines its orientation based on its Voronoi cell structure using a localized orientation approach. Since Voronoi vertices are equidistant from three or more sensors and are the candidates for the farthest point in a Voronoi cell, they naturally serve as candidate directions (see Theorem~\ref{t4}). To align each sensor for maximum coverage, the coverage area corresponding to its robust location is computed for each candidate vertex and the orientation that yields the highest average coverage is selected as the optimal direction.
To compute this area, we focus on the sensor's coverage within its own Voronoi cell. This results in four possible interaction cases between the sensor and its cell, as seen for sensors $s_1, s_2,s_3$ and $s_4$ in Figure~\ref{fig:geometric}. We analyze each of these cases individually.
\begin{figure}
\centering
\captionsetup{justification=centering}
\begin{subfigure}{.45\textwidth}
  \centering
  \includegraphics[width=.5\linewidth]{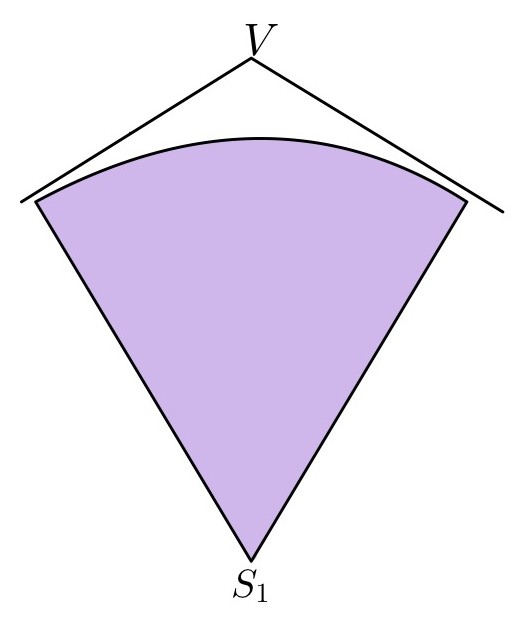}
  \caption{Complete sector inside Voronoi cell}
  \label{fig:localized1}
\end{subfigure}%
\begin{subfigure}{.45\textwidth}
  \centering
  \includegraphics[width=.5\linewidth]{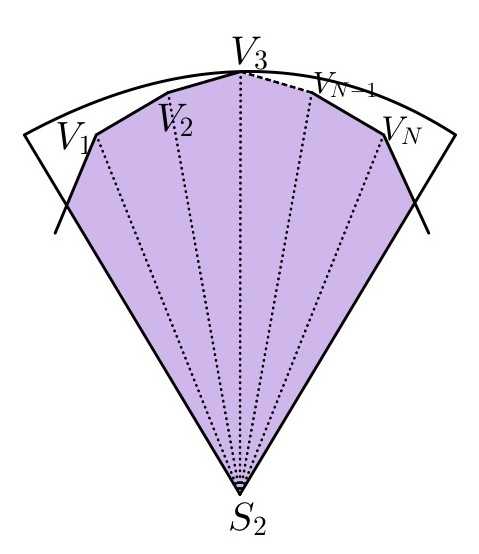}
  \caption{Partial Voronoi vertices inside sector}
  \label{fig:localized2}
\end{subfigure}%
\newline
\begin{subfigure}{.45\textwidth}
  \centering
  \includegraphics[width=.5\linewidth]{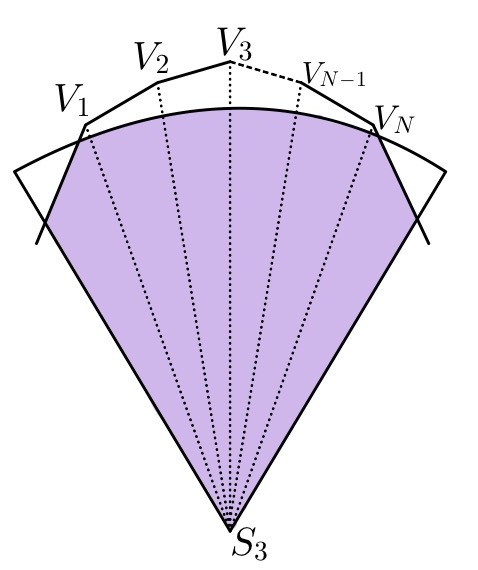}
  \caption{Voronoi edges intersect sector arc and sidelines}
  \label{fig:localized3}
\end{subfigure}%
\begin{subfigure}{.45\textwidth}
  \centering
  \includegraphics[width=.5\linewidth]{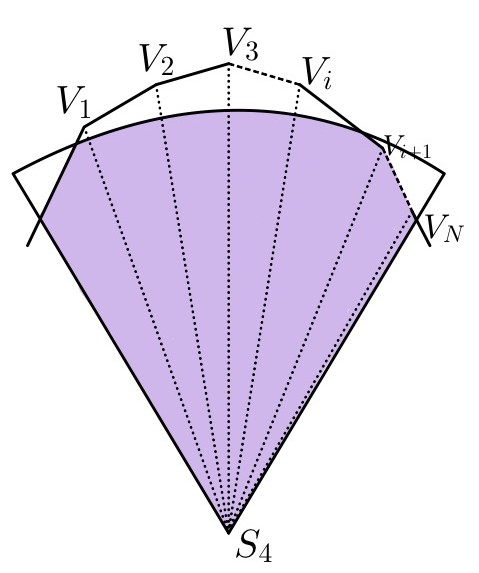}
  \caption{Mixed intersection (Cases 2 and 3 combined)}
  \label{fig:localized4}
\end{subfigure}
\caption{Illustration of the four possible cases of sensor coverage within a Voronoi cell}
\label{fig:voronoi_diag}
\end{figure}

{\bf Case-1: Complete sector inside Voronoi cell}\label{case1} [see Figure~\ref{fig:localized1}]. If the entire sensing sector of the $i^{th}$ sensor located at $s_i$ lies within its Voronoi cell, the coverage area $A_i$ is simply given by the sector area formula:
\begin{equation*}
    A_i= \frac{\theta_s}{2} r_s^2,
\end{equation*}
where $\theta_s$ is the sensor’s field of view in radians and $r_s$ is its effective coverage range. 

{\bf Case-2: Partial Voronoi vertices inside sector}\label{case2} [see Figure~\ref{fig:localized2}]. When some Voronoi vertices fall inside the sensing sector and two Voronoi edges intersect the sector’s sidelines, the coverage area is computed by dividing the region into $N+1$ triangles formed by connecting the sensor location to each vertex. The area of each triangle is obtained using Heron’s formula:
\begin{equation} \label{eqn2}
    \text{Area of a triangle}= \sqrt{d(d-e_1)(d-e_2)(d-e_3)},
\end{equation}
where $d= \frac{e_1+e_2+e_3}{2}$ is the semiperimeter and $e_j\;;\;j=1:3$ are the triangle’s edges.

To compute the required area, two types of edges are considered: (i) edges connecting two Voronoi vertices and (ii) edges connecting a Voronoi vertex with the intersection point of the sector sideline and a Voronoi edge. Given two vertices, say  $u = (u_1,u_2)$ and $v=(v_1,v_2)$, a point $w=(w_1,w_2)$ lying on the line segment joining them can be determined by:
\begin{equation} \label{eqn3}
    \begin{cases}
        w_2=v_2+\left(\frac{v_2-u_2}{v_1-u_1} \right) (w_1-v_1),\\
        d(u,w)+d(w,v)=d(u,v).
    \end{cases}
\end{equation}
The Voronoi vertices and edges are known from the Voronoi diagram, while the intersection points with the sector sidelines are derived by first formulating the sector sideline equations. Clearly, we get the sector sideline equation when the equality in equation~(\ref{cond2}) holds. Therefore, for a sensor nominally located at $s_i^o:=(x_i^o,y_i^o)$ with a robust position $s_i^w = (x_i^w,y_i^w)$ and the effective sensing range $r_i$, the sector sideline equations are given by:
\begin{equation} \label{eqn4}
\begin{cases}
    (a_o-x_i^w) \cos{o_i}+(b_0-y_i^w) \sin{o_i}= \sqrt{(a_0-x_i^w)^2 + (b_0-y_i^w)^2} \cos{\frac{\theta_s}{2}}, \\
    d(s_i^w,p) \leq r_i.
\end{cases}
\end{equation}

{\bf Case-3: Voronoi edges intersect sector arc and sidelines}\label{case3} [See Figure~\ref{fig:localized3}]. In this scenario, the area is computed in similar manner as Case 2 by constructing line segments between the sensor's robust location $s_i^w$ and the intersection points of the sector boundaries with the Voronoi edges. Given the sensor's robust position $s_i^w = (x_i^w,y_i^w)$, sensing range $r_i$, angle of view $\theta_s$ and orientation angle $o_i$, the area is calculated using equations~(\ref{eqn2}-\ref{eqn4}) along with the sector arc equation:
\begin{equation} \label{eqn5}
    \begin{cases}
        d(s_i^w,p) = r_i,\\
        a_0 \in [r_i \cos{\left( o_i + \frac{\theta_s}{2} \right)}+ x_i^w, r_i \cos{\left( o_i - \frac{\theta_s}{2} \right)} +x_i^w],\\
        b_0 \in [r_i \sin{\left( o_i - \frac{\theta_s}{2} \right)}+ y_i^w, r_i \sin{\left( o_i + \frac{\theta_s}{2} \right)} +y_i^w],
    \end{cases}
\end{equation}
where $p=(a_0,b_0)$ denotes a point located on the arc boundary of the sensor’s coverage sector.

{\bf Case-4: Mixed intersection} [See Figure~\ref{fig:localized4}]. In this case, the sensor's robust location is connected via line segments to all intersection points of the sector sensing region with Voronoi edges and the Voronoi vertices. The required area  $A_i$ is then determined using equations~(\ref{eqn2})-(\ref{eqn5}).\\
After computing the desired coverage area for all four cases, the sensor selects a unit direction vector pointing toward the vertex that provides the maximum coverage. This vector is chosen as the primary optimal orientation direction. Detailed steps are provided in Algorithm 1.

\begin{algorithm}
\caption{Localized Voronoi-Based Robust Orientation Optimization (LV-ROO)}
\label{algo:local_orientation}
\begin{algorithmic}[1]
\Require Sensor set $S$, Voronoi diagram $VD(S)$, uncertainty bound $\rho_i$, sensing radius $r_s$
\Ensure Optimal orientation $o_i^*$

\State Extract Voronoi vertices $\mathcal{V}_i$ and edges $\mathcal{E}_i$ of $VC(s_i)$

\State Define uncertainty-adjusted sensing region $\mathcal{B}_i(o_i)$ using RRF

\State Compute intersection set $\mathcal{I}_i = \partial \mathcal{B}_i(o_i) \cap VC(s_i)$

\State Discretize $VC(s_i)$ into subregions $\{\mathcal{F}_{ik}\}$ using $\mathcal{V}_i$, $\mathcal{E}_i$, and $\mathcal{I}_i$

\State Define candidate orientations $\Theta_i$ induced by $\mathcal{V}_i$ and $\{\mathcal{F}_{ik}\}$

\For{each $o_i \in \Theta_i$}
\State Compute covered region $\mathcal{C}_i(o_i) = \mathcal{B}_i(o_i) \cap V_i$
\State Evaluate area $A_i(o_i) = \sum_k |\mathcal{C}_i(o_i) \cap \mathcal{F}_{ik}|$
\EndFor

\State $o_i^* = \arg\max_{o_i \in \Theta_i} A_i(s_i,o_i)$

\State Update orientation $s_i \gets o_i^*$

\State Return optimized orientation for $s_i$
\end{algorithmic}
\end{algorithm}

\subsection{\bf{Collaborative Optimization and Boundary-Aware Refinement}}
\subsubsection{Collaborative Directional Adjustment}
\begin{figure}
\centering
\captionsetup{justification=centering}
\begin{subfigure}{.3\textwidth}
  \centering
  \includegraphics[width=.8\linewidth]{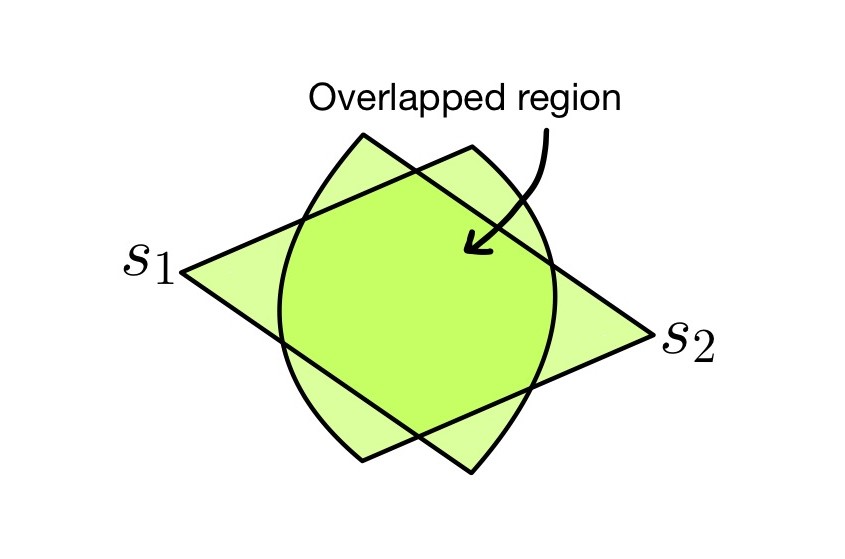}
  \caption{Sideline to sideline interaction}
  \label{fig:overlap1}
\end{subfigure}%
\begin{subfigure}{.3\textwidth}
  \centering
  \includegraphics[width=.8\linewidth]{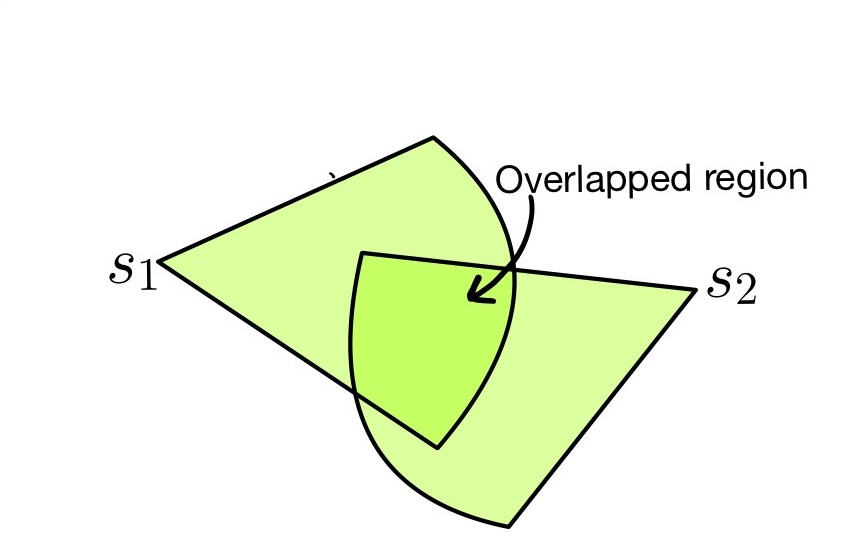}
  \caption{Sector arc to sideline interaction}
  \label{fig:overlap2}
\end{subfigure}%
\begin{subfigure}{.3\textwidth}
  \centering
  \includegraphics[width=.8\linewidth]{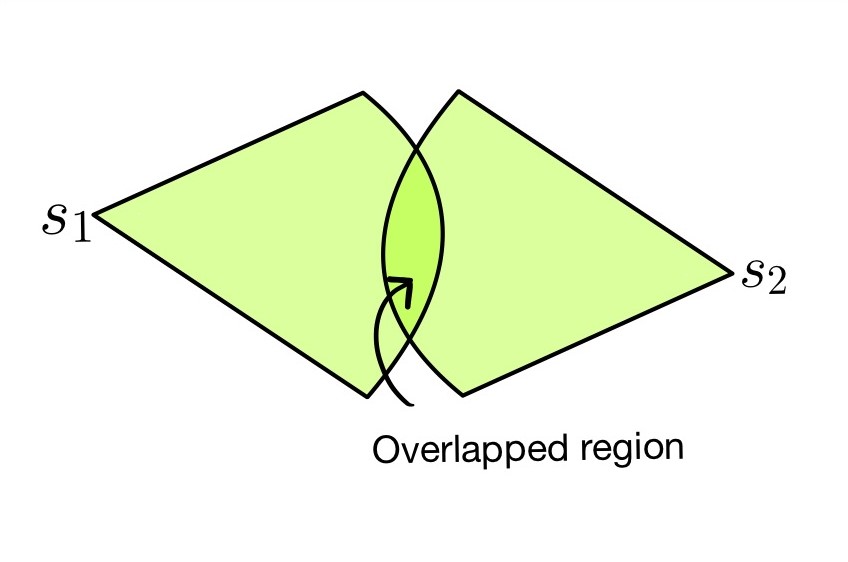}
  \caption{Sector arc to sector arc interaction}
  \label{fig:overlap3}
\end{subfigure}%
\caption{Possible interactions within different directional sensors}
\label{fig:collaborative1}
\end{figure}

After each sensor determines its localized optimal orientation, an inter-cell coordination step check for conflicts and minimizes sensing redundancy. This involves evaluating neighboring sensors' coverage within the same Voronoi partition and adjusting orientations based on the degree of overlap. The objective is to maximize coverage within each Voronoi cell while avoiding unnecessary redundancy. \\
For $i, j \in \left\{ 1,2 , \cdots, m\right\},$ consider two neighboring sensors at locations $s_i^o$ and $s_j^o$ with robust locations $s_i^w$ and $s_j^w$, robust feasibility radii $\rho_i$ and $\rho_j$, and effective sensing radii $r_i$ and $r_j$, respectively. Overlap is checked only if $d(s_i^w,s_j^w) \leq r_i+r_j$. For simplicity, we assume all sensors share the same sensing radius $r_s$, though the approach remains applicable to heterogeneous settings. If both sensors select the same vertex as their localized optimal orientation and $d(s_i^w,s_j^w) \leq r_i+r_j$, then overlap is inevitable.\\
When overlap is detected, one sensor reorients toward an uncovered region, prioritizing areas with higher coverage density. The sensor with smaller coverage switches to its next-best vertex until the conflict is resolved while checking for conflicts with other sensors. The resulting overlap can fall into one of three types, as demonstrated in Figure~\ref{fig:collaborative1}, with preferences ordered as \ref{fig:overlap1}, \ref{fig:overlap2} and \ref{fig:overlap3}, corresponding to decreasing levels of overlap. If all alternative vertices are exhausted during conflict resolution, the sensor adopts the orientation that yields minimal overlap.

\subsubsection{\bf{Perimeter Coverage Refinement}} 
Sensors near the network boundary refine their orientation to reduce outward coverage loss and improve interior sensing: if a sensor lies within a threshold $\epsilon$ of the boundary, it excludes vertices within $\epsilon$ of the region edges (vertical or horizontal). It removes Voronoi–boundary intersection points from its admissible orientation set, thereby promoting inward-facing directions and enhancing overall coverage efficiency. 

\subsection{\bf{Stopping Condition and Practical Considerations}} 
Since we aim to determine the worst-case scenario for real-world applications, we impose both upper and lower bounds on the RRF of each sensor. The RRF acts as a tightening factor on all coverage–feasibility constraints. For each candidate orientation, we compute the feasible sensing region under worst-case displacement within the RRF neighborhood, ensuring robustness by construction. In practical settings, the RRF of a sensor is finite, as beyond a certain distance from the nominal location, the sensor is guaranteed not to be deployed. Therefore, we define a maximum RRF value, $r_{max}$, beyond which further extension is unnecessary. Similarly, since the RRF must be above a certain threshold in real-life environments, we establish a minimum allowable value, $r_{min}$. As a constraint, if the computed RRF of a sensor exceeds $r_{max}$, it is set to $r_{max}$; if it falls below $r_{min}$, it is set to $r_{min}$. Since sensors cannot drift arbitrarily, the RRF is bounded: $$r_{min} \leq \rho _s \leq r_{max}.$$

Additional refinements can be incorporated into the algorithm. For instance, a stopping criterion can be introduced where the algorithm terminates once a predefined coverage requirement is met. Another approach could involve ceasing the addition of new sensors when the incremental coverage gain from deploying an additional sensor falls below a threshold $\delta$, ensuring efficiency. Mathematically, this condition can be expressed as
\begin{equation*}
    A_2 - A_1 < \delta,
\end{equation*}
where $A_1$ and $A_2$ represent the coverage area before and after placing the new sensor, respectively. The parameter $\delta$ regulates the trade-off between maximizing coverage and minimizing redundancy. However, such refinements are more relevant for structured deployment strategies rather than random sensor placement. Since our focus is solely on ensuring robust coverage under positional uncertainty, we do not incorporate these additional refinements in this study. Hence, if all vertex options are exhausted without resolving overlap during collaborative directional adjustment, the sensor is assigned the orientation with minimal overlap. Energy conservation via deactivation is not considered in this study.

\subsection{\bf{Algorithm for Coverage Optimization}}
We employ a Voronoi-based iterative optimization where each sensor orients toward the vertex that maximizes its worst-case coverage under uncertainty. Our approach iteratively refines sensor orientation and placement without requiring global information. The algorithm refines feasible vertices, optimizes local orientations, and resolves directional conflicts collaboratively. It includes an adaptive mechanism that dynamically adjusts sensor directions, countering deployment inaccuracies and enhancing overall coverage. The method operates using only local Voronoi information and follows three key principles: (i) each sensor optimizes its orientation based on local Voronoi structures to maximize immediate coverage, (ii) neighboring sensors collaboratively adjust their orientations to reduce coverage redundancy and fill sensing gaps, and (iii) sensors near the boundary refine their orientations to prevent excessive coverage spillover and enhance the useful sensing area.

The process begins by initializing the sensing region $R$, number of sensors $m$, and sensing angle $\theta_s$, followed by the construction of a Voronoi diagram to segment the sensing region based on nominal sensor placements. Using Theorem~\ref{theorem_main} and equation~(\ref{5.2}), sensors move to robust locations, after which orientations are optimized locally and updated iteratively until no two sensors share a vertex. If all reassignment options are exhausted, the sensor may retain its initial preferred orientation while still contributing to robustness. In the orientation selection step, instead of choosing the orientation corresponding to a single worst-case realization, we adopt an average-case robust criterion, wherein the selected orientation maximizes the average coverage evaluated over all worst-case realizations within the uncertainty set. This yields a less conservative and more balanced solution, improving practical coverage while still accounting for adverse positional deviations. The complete procedure is summarized in Algorithm~\ref{algo:coverage}.

\begin{algorithm}
\caption{Integrated Voronoi-Based Robust Orientation Optimization (IV-ROO)}
\label{algo:coverage}
\begin{algorithmic}[1]
\Require Sensor set $S=\{s_1,\ldots,s_m\}$ with orientations, uncertainty bounds $\rho$, sensing radius $r_s$
\Ensure Robust orientation maximizing coverage under positional uncertainty

\State Initialize sensor positions, orientations, and Voronoi regions; compute initial Voronoi diagram.

\State {\bf Step 1: Perimeter coverage refinement}
\For{each sensor $s_i \in S$}
\State Identify Voronoi vertices of $s_i$
\State Filter vertices that lie within boundary constraints $[\epsilon, R-\epsilon]$
\State Update valid vertex set for each sensor
\EndFor

\State {\bf Step 2: Localized orientation optimization}
\For{each sensor $s_i \in S$}
\State Compute the area of each possible orientation using Voronoi partitioning
\State Choose orientation maximizing the Voronoi cell area testing every orientation option under its RRF-adjusted footprint.
\State Update sensor orientation
\EndFor

\State {\bf Step 3: Collaborative directional adjustment}
\Repeat
\For{each pair of sensors $(s_i, s_j)$ where $s_i \text{ and } s_j$ share the same vertex}
\State Compute Euclidean distance $d_{ij}$ between $s_i$ and $s_j$
\If{$d_{ij} < \rho_i + \rho_j + 2r_s$}
\State Compute coverage area $A_i, A_j$
\State Sensor with lower coverage reorients to next best available vertex
\EndIf
\EndFor
\Until{All sensors have unique or exhausted orientations}

\State Compute final Voronoi diagram and compare coverage improvements
\State Generate performance metrics and analyze sensor placement effectiveness
\end{algorithmic}
\end{algorithm}

\section{Experimental Results and Analysis} 
This section evaluates the performance of the proposed framework under varying parameter conditions and uncertainty levels. We analyze the effectiveness of the robust orientation strategy through geometric visualization, quantitative coverage metrics, and resilience under sensor failures.

\subsection{\bf{Experimental Setup}} \label{6.2} 
All simulations are implemented in Python 3.12 using the Spyder environment on a standard laboratory workstation. Simulations are carried out on a domain of size $R$ of $(1000 \times 1000)$ square units. Sensors are deployed with predefined sensing radius $r_s$ and angle-of-view $\theta_s$, and their nominal positions are perturbed within bounded uncertainty sets characterized by RRF. Voronoi diagrams are constructed based on sensor locations to partition the sensing region, and all computations related to coverage and orientation optimization are performed within these partitions. Computational time for different deployment sizes is recorded using the built-in timing functions of Python. Experiments use $m=70$ sensors, $\rho\in[25,35]$, $\theta=60^{\circ}$, and $r_s=100$ unless mentioned otherwise. Performance is measured by analyzing the overall sensing coverage ratio under these conditions. To evaluate the effectiveness of our approach, we conduct simulations across various network configurations. 

Optimized sensor orientations under uncertainty are then obtained using Algorithm~\ref{algo:coverage}. We evaluate the coverage obtained by our proposed algorithms (LV-ROO and IV-ROO) and compare it with the initial coverage. We also present the performance upper bound with the "Oracle benchmark" coverage. This metric refers to an idealized case where exact future sensor positions are known a priori, allowing optimal orientation selection for those realizations. This serves as a theoretical benchmark, since such foresight is not available in practice, but helps quantify how close the proposed methods approach the optimal outcome. Simulations are done for two different sensor positions: exact case (nominal locations) and real-world case (perturbed locations). Each configuration is simulated for 500 iterations and the average results are reported to ensure a reliable and robust evaluation. 

\begin{figure}
\centering
\captionsetup{justification=centering}
\begin{subfigure}{.5\textwidth}
  \centering
  \includegraphics[width=.7\linewidth]{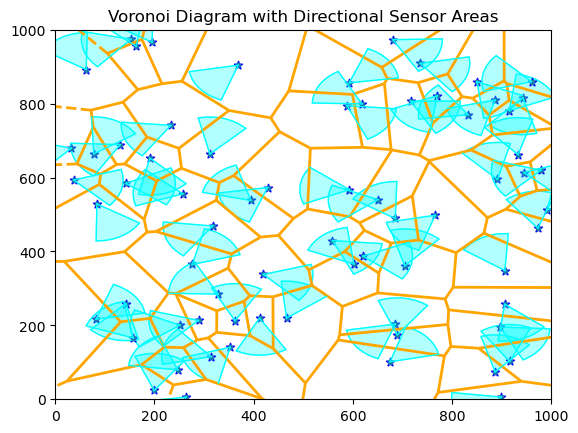}
  \caption{VD with initial orientation and nominal location}
  \label{fig:4.1a}
\end{subfigure}%
\hfill
\begin{subfigure}{.5\textwidth}
  \centering
  \includegraphics[width=.7\linewidth]{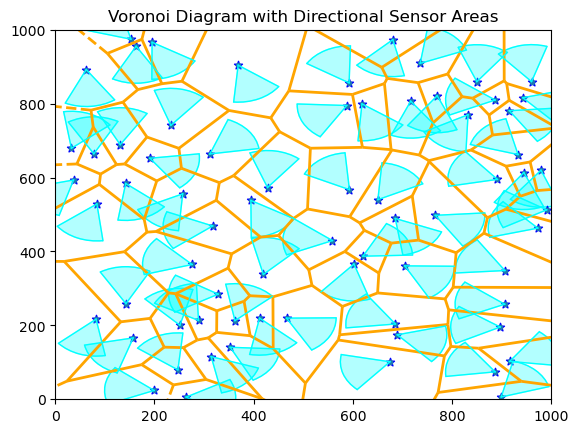}
  \caption{VD with robustified orientation and nominal location}
  \label{fig:4.1b}
\end{subfigure}

\vspace{0.5cm} 

\begin{subfigure}{.5\textwidth}
  \centering
  \includegraphics[width=.7\linewidth]{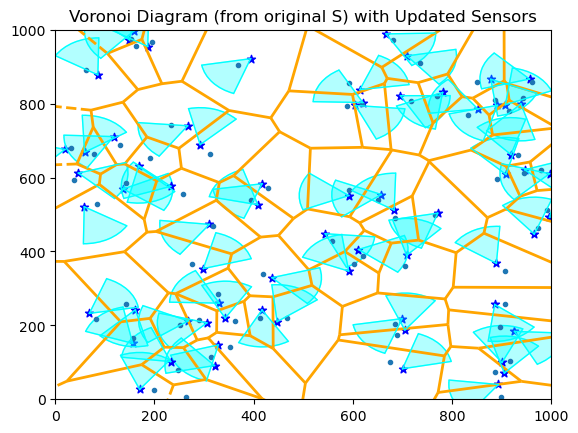}
  \caption{VD with initial orientation and perturbed location}
  \label{fig:4.1c}
\end{subfigure}%
\hfill
\begin{subfigure}{.5\textwidth}
  \centering
  \includegraphics[width=.7\linewidth]{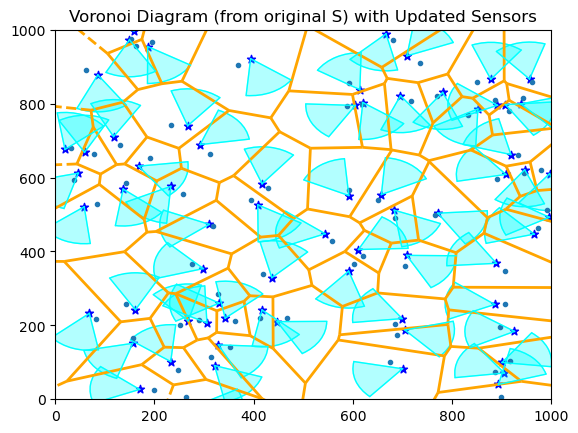}
  \caption{VD with robustified orientation and perturbed location}
  \label{fig:4.1d}
\end{subfigure}

\caption{Effect of Uncertainty on Sensor Orientation}
\label{VD_6.21}
\end{figure}

\subsection{\bf{Orientation and Coverage Behavior under Varying Conditions}}

We first analyze the impact of uncertainty on sensor orientation by comparing nominal and perturbed cases through Voronoi diagrams \ref{VD_6.21}. The comparison highlights the difference between initial (non-robust) orientations and robustified orientations obtained from the proposed method, demonstrating improved stability under positional variations.

\begin{figure}
\centering
\captionsetup{justification=centering}

\begin{subfigure}{.3\textwidth}
  \centering
  \includegraphics[width=.8\linewidth]{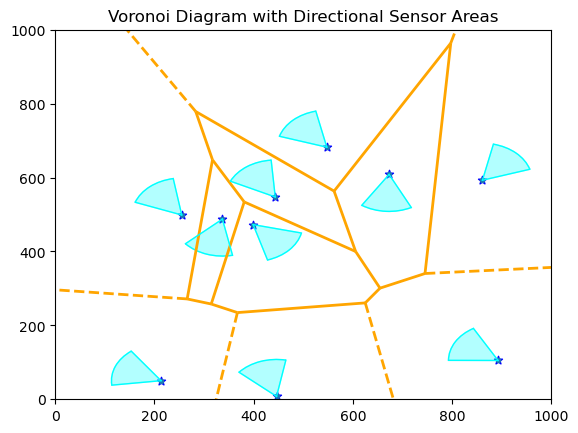}
  \caption{m = 10, t = 0.18 sec, \\coverage = 51161.47 sq units.}
\end{subfigure}%
\begin{subfigure}{.3\textwidth}
  \centering
  \includegraphics[width=.8\linewidth]{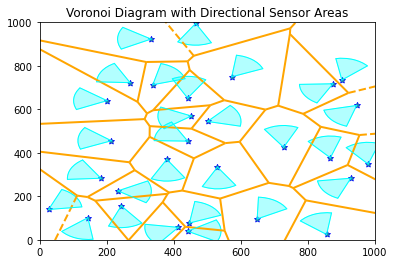}
  \caption{m = 30, t = 1.20 sec, \\coverage = 14599.39 sq units.}
\end{subfigure}%
\begin{subfigure}{.3\textwidth}
  \centering 
  \includegraphics[width=.8\linewidth]{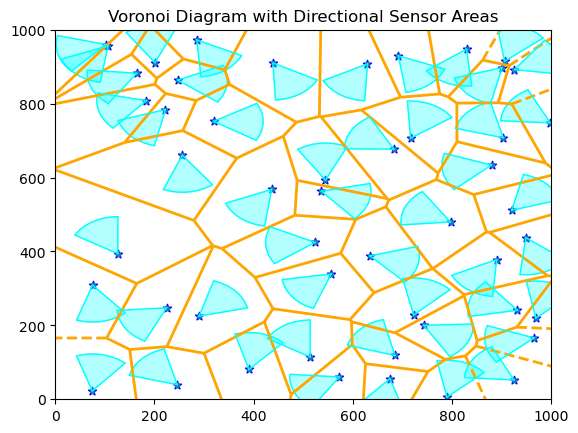}
  \caption{m = 50, t = 2.48 sec, \\coverage = 205520.81 sq units.}
\end{subfigure}%

\vspace{0.5cm} 

\begin{subfigure}{.3\textwidth}
  \centering
  \includegraphics[width=.8\linewidth]{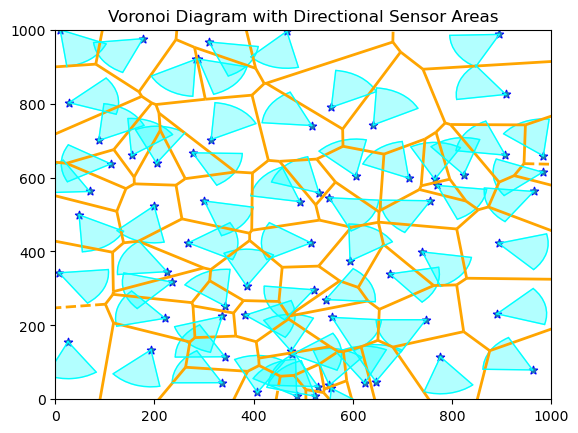}
  \caption{m = 70, t = 4.65 sec, \\coverage = 245831.66 sq units.}
\end{subfigure}%
\begin{subfigure}{.3\textwidth}
  \centering
  \includegraphics[width=.8\linewidth]{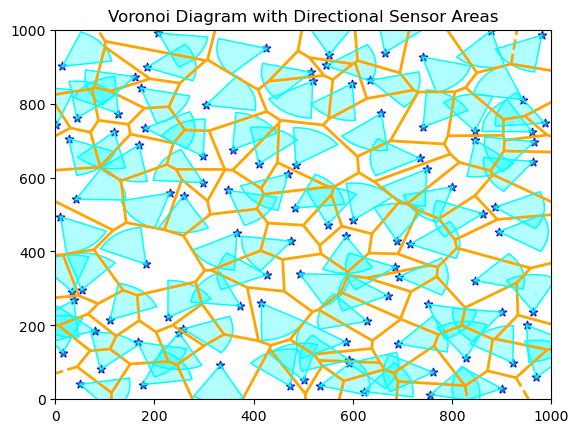}
  \caption{m = 100, t = 8.76 sec, \\coverage = 297641.67 sq units.}
\end{subfigure}%
\begin{subfigure}{.3\textwidth}
  \centering 
  \includegraphics[width=.8\linewidth]{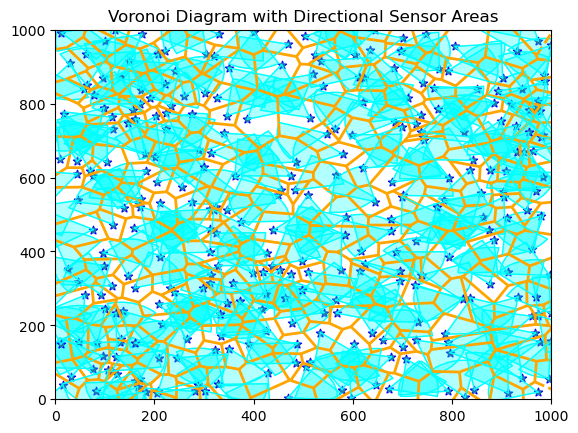}
  \caption{m = 200, t = 73.91 sec, \\coverage = 145993.39 sq units.}
\end{subfigure}%

\caption{Scalability Assessment with Varying Sensor Density}
\label{VD_6.22}
\end{figure}

{\bf {Impact of varying number of sensors ($m$)}:} Next, the scalability of the approach is examined by varying the number of deployed sensors ($m$); see Figure~\ref{VD_6.22}. For each case, Voronoi diagrams are generated to visualize spatial distribution and coverage patterns. As expected, increasing the number of sensors $m$ leads to higher coverage; however, beyond a certain density, the marginal gain becomes negligible due to increased overlap, indicating that further deployment offers limited practical benefit. The computational time required for each deployment size is also recorded, illustrating the efficiency of the proposed algorithm.

Further, a parametric analysis is conducted to study the influence of key parameters, including sensing angle ($\theta_s$), sensing radius ($r_s$), and uncertainty bounds (RRF). Coverage performance is evaluated for three configurations: initial (non-optimized), localized optimization (LV-ROO), and the integrated robust optimization (IV-ROO). The results are summarized in tabular form for systematic comparison.

{\bf {Impact of varying angle of view ($\theta$)}:}
Table~\ref{VD_6.231} reveals that as $\theta$ increases, coverage improves due to a wider sensing sector; however, the proposed LV-ROO and IV-ROO consistently achieve coverage closer to the ideal case under uncertainty, with IV-ROO showing better stability for smaller and moderate angles.

\begin{table}
\centering
\caption{Impact of varying angle of view $\theta_s$.}
\begin{tabular}{|p{2.2cm}|c|c|c|c|c|}
\hline
{\bf Varying $\theta_s$} & $\theta = 30^{\circ}$ & $\theta = 60^{\circ}$ & $\theta = 90^{\circ}$ & $\theta = 180^{\circ}$ & $\theta = 360^{\circ}$\\
\hline
\multicolumn{6}{|c|}{\bf Nominal Case} \\
\hline
{\bf Initial} & 73267.55 & 144464.51 & 217557.32 & 430443.01 & 861795.98 \\
{\bf LV-ROO}     & 134089.79 & 244504.80 & 334060.88 & 536067.04 & 861795.98 \\
{\bf IV-ROO}     & 148946.69 & 271748.49 & 371022.90 & 595956.94 & 861795.99 \\
\hline
\multicolumn{6}{|c|}{\bf Perturbed Case} \\
\hline
{\bf Initial} & 69973.73 & 142041.61 & 210803.74 & 422820.93 & 839460.81 \\
{\bf LV-ROO}     & 122302.76 & 226760.67 & 312621.64 & 511398.07 & 839460.81 \\
{\bf IV-ROO}     & 135927.51 & 251952.51 & 347384.29 & 567858.53 & 839460.99 \\
\hline
\multicolumn{6}{|c|}{\bf Oracle Benchmark} \\
\hline
{\bf Oracle Case}     & 151968.09 & 275604.14 & 375676.33 & 602141.82 & 861798.45 \\
\hline
\end{tabular}
\label{VD_6.231}
\end{table}

{\bf {Impact of varying sensing radius ($r_s$)}:}
Table~\ref{VD_6.232} reveals that the coverage increases with $r_s$ for all configurations, but under uncertainty, LV-ROO and especially IV-ROO maintain coverage significantly closer to the ideal case, demonstrating stronger robustness as sensing range grows.

\begin{table}
    \centering
    \caption{Impact of varying sensing radius ($r_s$ ) on deployment strategies.}
\begin{tabular}{|p{2.2cm}|c|c|c|c|c|}
\hline
{\bf Varying $r_s$} & $r_s=$ 60 &$r_s=$ 80  & $r_s=$ 100 & $r_s=$ 120 & $r_s=$ 140\\
\hline
\multicolumn{6}{|c|}{\bf Nominal Case} \\
\hline
{\bf Initial} & 88305.40 & 118900.44 & 144464.51 & 151107.40 & 164026.89 \\
{\bf LV-ROO}     & 111849.66 & 180973.83 & 244504.80 & 287614.77 & 310208.31 \\
{\bf IV-ROO}     & 124315.82 & 201191.45 & 271748.49 & 319455.96 & 344893.15 \\  
\hline
\multicolumn{6}{|c|}{\bf Perturbed Case} \\
\hline
{\bf Initial} & 81391.35 & 114040.44 & 142041.61 & 148735.74 & 163364.54 \\
{\bf LV-ROO}     & 105105.25 & 188579.76 & 226760.67 & 270363.61 & 295682.90 \\
{\bf IV-ROO}     & 116782.85 & 169763.69 & 251952.51 & 300442.79 & 328381.19 \\
\hline
\multicolumn{6}{|c|}{\bf Oracle Benchmark} \\
\hline
{\bf Oracle Case}     & 125795.90 & 203956.70 & 275604.14 & 323456.76 & 348839.74 \\
\hline
\end{tabular}
    \label{VD_6.232}
\end{table}

{\bf {Impact of varying radius of robust feasibility ($\rho$)}:}
Table~\ref{VD_6.233} reveals that as $\rho$ increases, coverage degrades due to larger positional uncertainty; however, LV-ROO and IV-ROO consistently preserve higher coverage compared to the initial configuration, remaining closer to the ideal achievable performance and ensuring reliable sensing.

\begin{table}
    \centering
    \caption{Impact of varying radius of robust feasibility on deployment strategies.}
\begin{tabular}{|p{2.2cm}|c|c|c|c|c|}
\hline
{\bf Varying RRF} & 5 $\leq \rho \leq$ 15& 15$\leq \rho \leq$25 & 25$\leq \rho \leq$ 35& 35$\leq \rho \leq$ 45& 45$\leq \rho \leq$ 55\\
\hline
\multicolumn{6}{|c|}{\bf Nominal Case} \\
\hline
{\bf Initial} & 142189.00 & 142384.36 & 144464.51 & 239350.29 & 144948.69 \\
{\bf LV-ROO}     & 246522.08 & 245911.82 & 244504.80 & 266004.54 & 232236.25 \\
{\bf IV-ROO}     & 274069.90 & 273252.13 & 271748.49 & 274497.66 & 257895.40 \\  
\hline
\multicolumn{6}{|c|}{\bf Perturbed Case} \\
\hline
{\bf Initial} & 142855.25 & 141616.24 & 142041.61 & 133636.40 & 132475.84 \\
{\bf LV-ROO}     & 246374.48 & 241097.80 & 226760.67 & 214886.70 & 196182.78 \\
{\bf IV-ROO}     & 273589.64 & 267863.99 & 251952.51 & 238618.75 & 217945.22 \\
\hline
\multicolumn{6}{|c|}{\bf Oracle Benchmark} \\
\hline
{\bf Oracle Case}     & 274129.39 & 273968.38 & 275604.14 & 274497.66 & 274071.34 \\
\hline
\end{tabular}
    \label{VD_6.233}
\end{table}

\subsection{\bf{Comparative Coverage Analysis}}
\begin{figure}[h!]
    \centering
    \begin{subfigure}{0.48\textwidth}
        \centering
        \includegraphics[width=\linewidth]{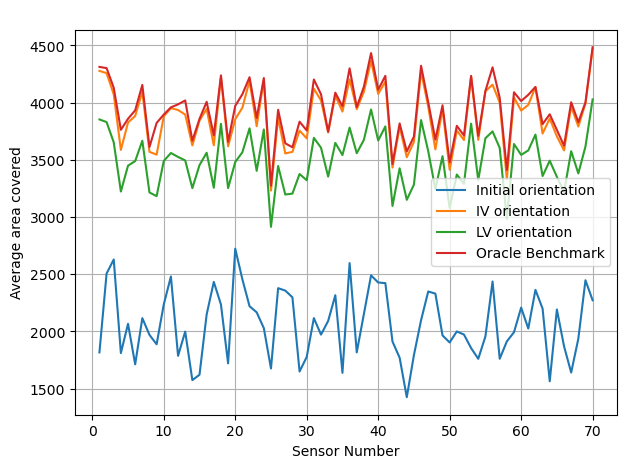}
        \caption{Coverage comparison under nominal orientations.}
        \label{fig:nominal_comparison}
    \end{subfigure}
    \hfill
    \begin{subfigure}{0.48\textwidth}
        \centering
        \includegraphics[width=\linewidth]{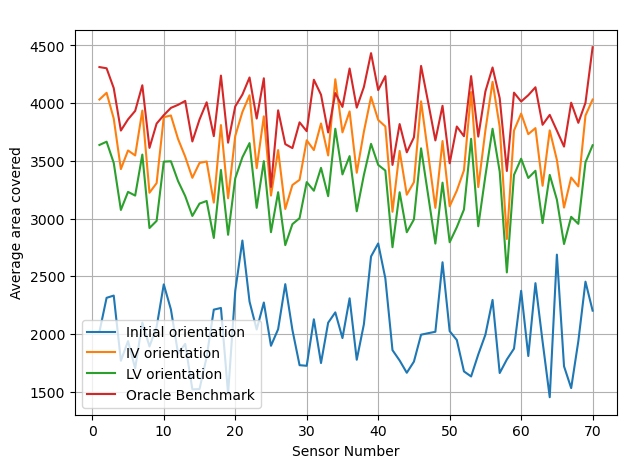}
        \caption{Coverage comparison under perturbed orientations.}
        \label{fig:worst_comparison}
    \end{subfigure}
    \caption{Voronoi-cell wise area distribution across sensors for different orientation strategies}
    \label{6.31}
\end{figure}

A detailed comparison of coverage performance is carried out across different configurations. First, a Voronoi cell-wise analysis is presented to quantify the contribution of each sensor under different strategies; see Figure~\ref{6.31}. This highlights how local and collaborative optimization improve coverage distribution.

\begin{figure}
\centering
\begin{subfigure}{0.48\columnwidth}
    \centering
    \includegraphics[width=\linewidth]{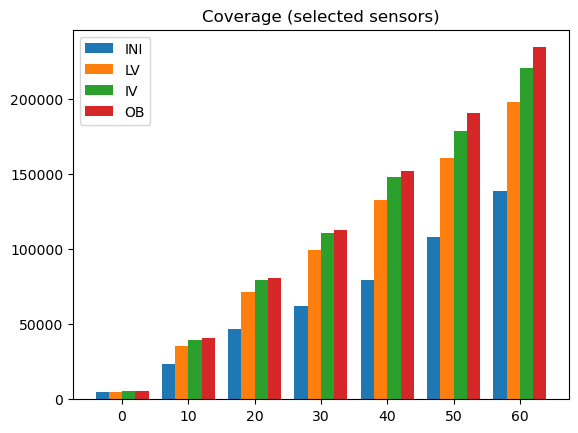}
    \caption{Scalability Performance Analysis}
    \label{6.32}
\end{subfigure}
\hfill
\begin{subfigure}{0.48\columnwidth}
    \centering
    \includegraphics[width=\linewidth]{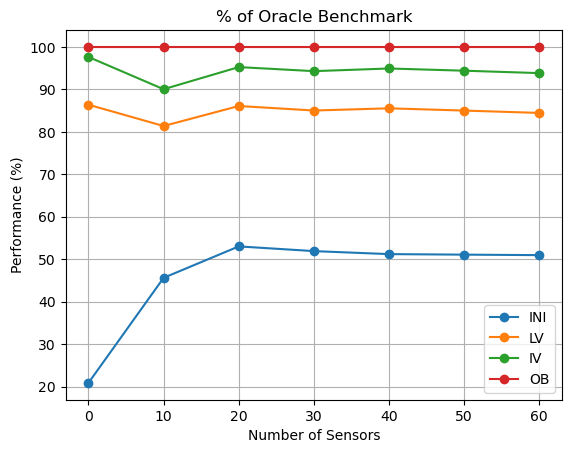}
    \caption{Relative Performance Comparison}
    \label{6.33}
\end{subfigure}
\caption{Performance Evaluation}
\label{fig:coverage_graphs}
\end{figure}

Next, the total coverage is plotted against the number of sensors, demonstrating the scalability and efficiency of each approach as shown in Figure~\ref{6.32}. In addition, performance is evaluated in terms of percentage coverage relative to the "Oracle benchmark" coverage, see Figure~\ref{6.33}, computed as the ratio of achieved coverage to the ideal "Oracle benchmark" coverage under known perturbed positions, i.e., $$\textbf{Performance metric}= \frac{\textbf{Method Coverage}}{\textbf{"Oracle benchmark Coverage"}} * 100$$ Bar charts are used to compare all configurations, demonstrating the effectiveness of the proposed robust optimization framework in improving both absolute and relative coverage.

\subsection{\bf{Robustness under Sensor Failures}}
\begin{figure}
    \centering
    \includegraphics[width=0.5\textwidth]{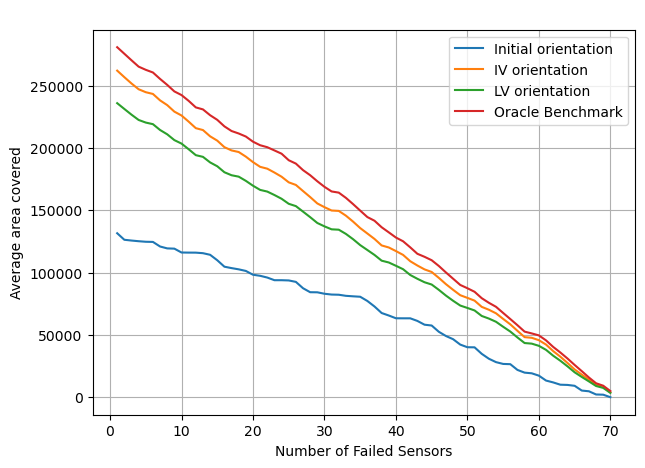}
    \caption{Resilience Under Sensor Failures}
    \label{6.4}
\end{figure}

To evaluate resilience, we analyze the impact of sensor failures on overall coverage. Different scenarios are considered where a subset of sensors becomes non-functional. Coverage degradation is studied for all configurations, along with an additional ideal case representing the "Oracle benchmark" performance as demonstrated in Figure~\ref{6.4}. Graphs are presented to show the relationship between the number of failed sensors and the resulting coverage. Furthermore, comparative plots illustrate how each method performs relative to the ideal case as failures increase. The results demonstrate that the proposed robust framework maintains higher coverage levels and exhibits graceful degradation under adverse conditions.

\section{Conclusions and Future Scopes}
In this work, we developed a robust Voronoi-based framework for sensor orientation optimization in directional sensor networks under positional uncertainty, where the RRF is integrated to explicitly quantify and handle deployment perturbations. The proposed approach combines localized orientation optimization, collaborative conflict resolution, and boundary-aware refinement within a unified structure that relies solely on local Voronoi information, ensuring both scalability and computational efficiency. The results demonstrate that the framework consistently enhances robustness and coverage stability, while enabling cooperative orientation adjustments among sensors to maintain reliable coverage under deviations. Although a slight reduction in nominal coverage is observed, the method significantly outperforms initial configurations and baseline approaches in uncertain settings, exhibiting graceful performance degradation under sensor failures and improved resilience to environmental variability.

Overall, the proposed framework provides a practically viable and theoretically grounded solution that bridges the gap between robustness modeling and real-world coverage optimization in directional sensor networks. By leveraging RRF as a principled robustness measure, the approach ensures tolerance to misplacements and delivers stable performance in non-ideal deployment scenarios, making it well-suited for real-world monitoring applications. Future work can extend this framework to heterogeneous sensor models, incorporate learning-based predictive orientation strategies, and adapt the formulation to multi-terrain environments through weighted distance metrics and region-aware feasibility constraints, thereby further improving scalability, adaptability, and real-world applicability.

\end{document}